\documentclass[12pt]{article}  
\def\sq{\hbox {\rlap{$\sqcap$}$\sqcup$}}
\def\sq{\hbox {\rlap{$\sqcap$}$\sqcup$}}
\def\R{ {\rm R \kern -.31cm I \kern .15cm}}
\def\C{ {\rm C \kern -.15cm \vrule width.5pt \kern .12cm}}
\def\Z{ {\rm Z \kern -.27cm \angle \kern .02cm}}
\def\N{ {\rm N \kern -.26cm \vrule width.4pt \kern .10cm}}
\def\1{{\rm 1\mskip-4.5mu l} }
\def\lsim{\raise0.3ex\hbox{$<$\kern-0.75em\raise-1.1ex\hbox{$\sim$}}}
\def\gsim{\raise0.3ex\hbox{$>$\kern-0.75em\raise-1.1ex\hbox{$\sim$}}}
\def\noi{\noindent}

\def\beq{\begin{equation}}   \def\eeq{\end{equation}}
\def\bea{\begin{eqnarray}}  \def\eea{\end{eqnarray}}
\def\nn{\nonumber}
\def\noi{\noindent}
\def\beeq{\begin{eqnarray}} \def\eeeq{\end{eqnarray}}
\newcommand\mysection{\setcounter{equation}{0}\section}

\newcounter{hran}

\begin{document} 
\centerline{\large\bf Modified wave operators without loss of regularity} 
 \vskip 3 truemm 
 \centerline{\large\bf  for some long range Hartree equations. I} 
  \vskip 0.8 truecm

\centerline{\bf J. Ginibre}
\centerline{Laboratoire de Physique Th\'eorique\footnote{Unit\'e Mixte de
Recherche (CNRS) UMR 8627}}  \centerline{Universit\'e de Paris XI, B\^atiment
210, F-91405 ORSAY Cedex, France}
\vskip 3 truemm

\centerline{\bf G. Velo}
\centerline{Dipartimento di Fisica, Universit\`a di Bologna}  \centerline{and INFN, Sezione di
Bologna, Italy}

\vskip 1 truecm

\begin{abstract}
We reconsider the theory of scattering for some long range Hartree equations with potential $|x|^{-\gamma}$ with $1/2 < \gamma < 1$. More precisely we study the local Cauchy problem with infinite initial time, which is the main step in the construction of the modified wave operators. We solve that problem in the whole subcritical range without loss of regularity between the asymptotic state and the solution, thereby recovering a result of Nakanishi. Our method starts from a different parametrization of the solutions, already used in our previous papers. This reduces the proofs to energy estimates and avoids delicate phase estimates.
\end{abstract}

\vskip 1 truecm
\noi 2000 MSC :  Primary 35P25. Secondary 35B40, 35Q40, 81U99.\par \vskip 2 truemm

\noi Key words : Long range scattering, wave operators, Hartree equation. \par 
\vskip 1 truecm

\noindent LPT Orsay 12-33\par
\noindent April 2012\par \vskip 3 truemm

\newpage
\pagestyle{plain}
\baselineskip 18pt
\mysection{Introduction}
\hspace*{\parindent} 
This paper is devoted to the theory of scattering and more precisely to the proof of existence of modified wave operators for the long range Hartree type equation 
\beq
i \partial_t u = - (1/2) \Delta u + g(u) u 
\label{1.1e}
\eeq

\noi where $u$ is a complex valued function defined in space time ${I\hskip-1truemm R}^{n+1}$ with $n \geq 2$, 
$\Delta$ is the Laplace operator in ${I\hskip-1truemm R}^n$ and 
\beq
\label{1.2e}
g(u) = \kappa \  |x|^{-\gamma} \star |u|^2
\eeq

\noi where $\kappa \in {I\hskip-1truemm R}$, $0 < \gamma \leq 1$ and $\star$ denotes the convolution in ${I\hskip-1truemm R}^n$. \par

A large amount of work has been devoted to the theory of scattering for the Hartree equation (\ref{1.1e}) with nonlinearity  (\ref{1.2e}), both in the short range case $\gamma > 1$ and in the long range case $\gamma \leq 1$. See \cite{2r} \cite{3r} \cite{5r} \cite{6r} and references therein quoted. In order to prove the existence of wave operators, one has to construct solutions of the given equation with prescribed asymptotic behaviour at $\pm \infty$ in time. The asymptotic behaviour is that of solutions of the free Schr\"odinger equation in the short range case $\gamma > 1$, thereby leading to ordinary wave operators, and has to be modified by a suitable phase factor in the long range case $\gamma \leq 1$, thereby leading to modified wave operators in that case. The asymptotic behaviour is parametrized by an asymptotic state, which plays the role of (in the short range case can be taken to be) the initial data for the asymptotic behaviour. The main step in the construction of solutions with prescribed asymptotic behaviour consists in solving the local Cauchy problem with infinite initial time, with the asymptotic state playing the role of initial data. From now on, we concentrate on that problem. In \cite{2r} we have solved that problem for arbitrarily large data in the range $1/2 < \gamma < 1$ (the easier borderline case $\gamma = 1$ can be treated by the same method but requires slightly different formulas containing logarithms). The method used in \cite{2r} consists in parametrizing the solution $u$ in terms of an amplitude $v$ and a phase $\varphi$ and solving an auxiliary system of equations for the pair $(v, \varphi )$. It has two defects, namely~: (i) there is a natural notion of criticality for that problem, and the regularity required for the solution is significantly higher than the critical one, namely the method does not cover the entire subcritical range, and (ii) there occurs a loss of regularity (essentially a loss of two derivatives on $v$) between the asymptotic state and the solution eventually obtained. These two defects were remedied in \cite{5r} and \cite{6r} for $\gamma = 1$ and for $1/2 < \gamma < 1$ respectively, by the use of two new ingredients, namely~: (i) a different parametrization of the solution $u$, and (ii) the clever use of a local conservation law satisfied by Schr\" odinger type equations including (\ref{1.1e}). This allows in addition to fix the new phase $\varphi$ from the beginning, thereby leaving only one equation for the new amplitude $v$. However the method requires rather complicated phase estimates in the more difficult case $1/2 < \gamma < 1$. \par

It turns out however that the two new ingredients used in \cite{5r} \cite{6r} are independent of each other. In particular the local conservation law can also be exploited with the parametrization used in \cite{2r}, provided the latter is supplemented with the systematic use of an ultraviolet momentum cut-off. One can then recover the main results of \cite{6r}, namely solve the local Cauchy problem at infinity in time in the whole subcritical range and without any regularity loss, by elementary energy methods not requiring any delicate phase estimates. The purpose of the present paper is to present that simpler alternative method.\par

The simpler method also lends itself naturally to an iteration scheme which can be expected to cover the range $\gamma \leq 1/2$, with the $n$-th interation covering the range $1/(n+2) < \gamma < 1/(n+1)$, still without any regularity loss between the asymptotic state and the solution. However the method does not seem to make it possible to cover the entire subcritical range as soon as $\gamma \leq 1/2$, and stronger regularity conditions seem to be required. Furthermore, the necessary estimates, although still elementary, become more and more cumbersome as $n$ increases. In a subsequent paper, as an illustration we shall treat the problem in the range $1/3 < \gamma < 1/2$ by means of the first iteration.\par

We now introduce the relevant parametrization of $u$ needed to study the Cauchy problem at infinite time, restricting our attention to positive time. The unitary group 
\beq
\label{1.3e}
U(t) = \exp (i(t/2)\Delta ) \eeq

\noi which solves the free Schr\"odinger equation can be written as 
\beq
U(t) = M(t) \ D(t) \ F\ M(t) 
\label{1.4e}
\eeq

\noi where $M(t)$ is the operator of multiplication by the function
\beq
M(t) = \exp (ix^2/2t) \ , 
\label{1.5e}
\eeq

\noi $F$ is the Fourier transform and $D(t)$ is the dilation operator
\beq
D(t) = (it)^{-n/2} \ D_0(t) 
\label{1.6e}
\eeq

\noi where
\beq
\left ( D_0 (t) f \right ) (x) = f(x/t)\ . 
\label{1.7e}
\eeq

\noi For any function $w$ of space time, we define
\beq
\label{1.8e}
\widetilde{w}(t) = U(-t) \ w(t) 
\eeq

\noi and we define the pseudoconformal inverse $w_c$ of $w$ by
\beq
w(t) = M(t) \ D(t) \ \overline{w_c} (1/t) 
\label{1.9e}
\eeq

\noi or equivalently
\beq
\widetilde{w}(t) = \overline{F \widetilde{w}_c (1/t)}\ 
\label{1.10e}
 \eeq

\noi which shows that the pseudoconformal inversion is involutive. \par

The Cauchy problem at infinite initial time for $u$ is then equivalent to the Cauchy problem at initial time zero for its pseudoconformal inverse $u_c$. The equation (\ref{1.1e}) is replaced by 
\beq
i \partial_t  u_c = - (1/2) \Delta u_c + t^{\gamma -2} \ g(u_c) u_c \ .
\label{1.11e}
 \eeq

\noi We now parametrize $u_c$ in terms of an amplitude $v$ and a phase $\varphi$ according to 
\beq
\label{1.12e}
u_c(t) = \exp (- i \varphi (t)) v(t) 
\eeq

\noi so that
\begin{eqnarray*}
u(t) &=& M(t)\ D(t) \exp \left ( i \varphi (1/t)\right ) \overline{v}(1/t) \\
&=& D(t) \exp \left ( i \varphi (1/t)\right ) D^{-1}(t)\ M(t)\ D(t) \ \overline{v}(1/t)
\end{eqnarray*}

\noi or equivalently
\beq
\label{1.13e}
u(t) = \exp \left ( i \left ( D_0(t) \varphi (1/t)\right ) \right ) v_c(t)\ .
\eeq

\noi The original equation then becomes the following equation for $v$
\beq
\label{1.14e}
i \partial_t v = - (1/2) \Delta_s v + \left ( t^{\gamma -2} g(v) - \partial_t \varphi \right ) v 
\eeq

\noi where $s = \nabla \varphi$ and
\beq
\Delta_s = (\nabla - i\ s)^2 = \Delta - 2 i\ s \cdot \nabla - i (\nabla \cdot s) - |s|^2 \ .
\label{1.15e}
\eeq

We want to choose $\varphi$ so as to cancel the divergence at $t = 0$ of the last term in (\ref{1.14e}), but that cancellation is needed only at large distances, namely for low momentum. We therefore introduce a momentum cut-off as follows. Let $\chi \in {\cal C}^\infty ({I\hskip-1 truemm R}^+,{I\hskip-1 truemm R}^+)$, $0 \leq \chi \leq 1$, $\chi (\ell ) = 1$ for $\ell \leq 1$, $\chi (\ell ) = 0$ for $\ell \geq 2$. We define 
\beq
\label{1.16e}
\chi_L = \chi (\omega t^{1/2}) \quad , \quad \chi_S = 1 - \chi_L
\eeq

\noi with $\omega = (-\Delta )^{1/2}$, and correspondingly 
\beq
\label{1.17e}
g_L(v) = \chi_L \ g(v) \quad , \quad g_S(v) = \chi_S \ g(v) \ .
\eeq

\noi We want to solve (\ref{1.14e}) with $v$ continuous at $t=0$ with $v(0) = v_0$. For that purpose, we choose (assuming from now on $\gamma < 1$)
\beq
\label{1.18e}
\varphi = - (1 - \gamma )^{-1} \ t^{\gamma - 1} \ g_L (v_0)
\eeq

\noi so that 
\beq
\label{1.19e}
\partial_t \varphi = t^{\gamma - 2} \ g_L (v_0) - (1 - \gamma )^{{-1}}\ t^{\gamma - 2} \ \widetilde{\chi}_L \ g(v_0)
\eeq

\noi with
\beq
\label{1.20e}
\widetilde{\chi}_L = \widetilde{\chi} (\omega t^{1/2}) \ ,\widetilde{\chi} (\ell ) = (1/2) \ell \chi '(\ell ) \ .
\eeq

\noi With that choice, the equation (\ref{1.14e}) becomes
\beq
\label{1.21e}
i \partial_t v = L(v) \ v  
\eeq

\noi with
\beq
\label{1.22e}
L(v) = - (1/2) \Delta_s + t^{\gamma -2} \ g_S (v) + (1- \gamma )^{-1} \ t^{\gamma - 2} \  \widetilde{\chi}_L\ g(v_0) + t^{\gamma - 2} \left ( g_L (v) - g_L (v_0)\right ) \ .
\eeq

\noi We shall also need the partly linearized equation for $v'$
\beq
\label{1.23e}
i \partial_t v' = L(v) v'\ .  
\eeq

The method consists in first solving the Cauchy problem with initial time zero for the linearized equation (\ref{1.23e}). One then shows that the map $v \to v'$ thereby defined is a contraction in a suitable space in a sufficiently small time interval. This solves the Cauchy problem with initial time zero for the nonlinear equation (\ref{1.21e}). One then translates the results through the change of variables (\ref{1.12e}) to solve the Cauchy problem with initial time zero for the equation (\ref{1.11e}) or equivalently with infinite initial time for the equation (\ref{1.1e}). The final result can be stated as the following proposition, which is adapted to the equation (\ref{1.1e}) in a neighborhood of infinity in time. We need the notation
$$FH^\rho = \{ u \in {\cal S} ' : F^{-1} u \in H^\rho \} \ .$$

\noi {\bf Proposition 1.1.} {\it Let $1/2 < \gamma < 1$. Let $1 - \gamma /2 < \rho < n/2$. \par

(1) Let $u_0 \in F H^\rho$ and define
\beq
\label{1.24e}
\varphi (t) = - (1 - \gamma )^{-1} \ t^{\gamma -1} \ g_L (Fu_0) \ .
\eeq

\noi Then there exists $T_\infty > 0$ and there exists a unique solution $u$ of the equation (\ref{1.1e}) such that $v_c$ defined by (\ref{1.9e}) (\ref{1.12e}) or equivalently by (\ref{1.13e}) satisfies $\widetilde{v}_c \in {\cal C} ([T_\infty , \infty ), FH^\rho )$ and such that
\beq
\label{1.25e}
\widetilde{v}_c (t) \to u_0 \hbox{ in } FH^\rho \hbox{ when } t \to \infty \ .
\eeq

\noi Furthermore $\widetilde{u} \in {\cal C} ([T_\infty , \infty ), FH^\rho )$ and $\widetilde{u}$ satisfies the estimate 
\beq
\label{1.26e}
\parallel \widetilde{u} (t) ; FH^{\rho} \parallel \ \leq C\ a_0 \left ( 1 + a_0^2 \ t^{1 - \gamma} \right )^{1 + [\rho ]} 
\eeq

\noi for all $t \geq T_\infty$, where $[\rho ]$ is the integral part of $\rho$ and 
$$a_0 = \ \parallel u_0 ; FH^\rho \parallel \ .$$

(2) Let in addition $\rho > 3/4$. Then the map $u_0 \to \widetilde{v}_c$ is continuous from $FH^\rho$ to $L^\infty ([T_\infty , \infty ), FH^\rho )$ and the map $u_0 \to \widetilde{u}$ is continuous from $FH^\rho$ to $L^\infty ([T_\infty , T_1 ], FH^\rho )$ for all $T_1$, $T_\infty < T_1 < \infty$.}\\

Proposition 1.1 follows from Propositions 4.1, 4.2, 5.1 and 5.2 through the change of variables (\ref{1.9e}) or (\ref{1.10e}), which implies in particular that
\beq
\label{1.27e}
\parallel \widetilde{w} (t) ; FH^\rho \parallel \ = \ \parallel \widetilde{w}_c (1/t) ; H^\rho \parallel \ = \ \parallel w_c (1/t) ; H^\rho \parallel \ .
\eeq

\noi The condition $\rho > 1 - \gamma /2$ in Proposition 1.1 is the subcriticality condition mentioned above. For technical reasons, the continuity properties in Part (2) are proved only under the stronger condition $\rho > 3/4$. \par

This paper is organized as follows. In Section~2, we introduce some notation and we collect a number of estimates which are used throughout this paper. In Section~3, we study the Cauchy problem for the linearized equation (\ref{1.23e}) with initial time $t_0 \geq 0$. In Section~4, we solve the Cauchy problem with initial time zero for the nonlinear equation (\ref{1.21e}). In Section~5, we prove the continuity of the solutions of (\ref{1.21e}) with respect to the initial data.  In Appendix A1, we prove Lemma 2.4 on commutators. In Appendix A2, we sketch the proof of the preliminary Proposition 3.1.\par

In all this paper (as in \cite{2r}) we assume that $1/2 < \gamma < 1$. The easier case $\gamma = 1$ can be treated by the same method, but requires slightly different formulas.

\mysection{Notation and preliminary estimates} 
\hspace*{\parindent}
In this section we introduce some notation and we collect a number of estimates which will be used throughout this paper. We denote by
$\parallel \cdot \parallel_r$ the norm in $L^r \equiv L^r({I\hskip-1truemm R}^{n})$.  For any interval $I$ and any Banach space $X$ we denote by ${\cal C}(I,X)$ (resp. ${\cal C}_w(I,X))$ the space of strongly (resp. weakly) continuous functions from $I$ to $X$ and by $L^{\infty} (I, X)$ the space of measurable essentially bounded functions from $I$ to $X$. For real numbers $a$ and $b$ we use the notation $a \vee b = {\rm Max}(a,b)$ and $a\wedge b = {\rm Min} (a,b)$. We define $(a)_+ = a \vee 0$ and 
$$\begin{array}{lll}
[a]_+ &= (a)_+ &\qquad \hbox{for $a\not= 0$} \\
&& \\
&= \varepsilon \hbox{ for some $\varepsilon > 0$} &\qquad \hbox{for $a= 0$} \ .
\end{array}$$
\par

We shall use the Sobolev spaces $\dot{H}_r^\sigma$ and $H_r^\sigma$ defined for $- \infty < \sigma < + \infty$, $1 \leq r < \infty$ by
$$\dot{H}_r^\sigma = \left \{ u:\parallel u;\dot{H}_r^\sigma\parallel \ \equiv \ \parallel \omega^\sigma u\parallel_r \ <
\infty \right \}$$

\noi and

$$H_r^\sigma = \left \{ u:\parallel u;H_r^\sigma\parallel \ \equiv \ \parallel <\omega>^\sigma u\parallel_r \ <
\infty \right \}$$

\noi where $\omega = (- \Delta)^{1/2}$ and $< \cdot > = (1 + |\cdot |^2)^{1/2}$. The subscript $r$ will
be omitted both in $H^\sigma$ and in the $L^r$ norm if $r = 2$ and we shall use
the notation 
$$\parallel \omega^{\sigma \pm 0} u \parallel \ = \left ( \parallel \omega^{\sigma + \varepsilon} u \parallel\ \parallel \omega^{\sigma - \varepsilon} u \parallel\right )^{1/2} \quad \hbox{for some $\varepsilon > 0$\ .}$$  

Note also that for $0 < \gamma < n$ \cite{7r} 
$$g(u) = \kappa \ |x|^{-\gamma} \star |u|^2 = \kappa \ C_{\gamma , n} \ \omega^{\gamma - n} \ |u|^2 \ .$$ 

We shall use extensively the following Sobolev inequalities. \\

\noi {\bf Lemma 2.1.} {\it Let $1 < q, r < \infty$, $1 < p \leq \infty$ and $0 \leq \sigma < \rho$. If $p = \infty$, assume that $\rho - \sigma > n/r$. Let $\theta$ satisfy $\sigma /\rho \leq \theta \leq 1$ and
$$n/p - \sigma = (1 - \theta )n/q + \theta (n/r - \rho ) \ .$$

\noi Then the following inequality holds}
\beq  
\label{2.1e}
\parallel \omega^\sigma u\parallel_p \ \leq C  \parallel u \parallel_q^{1 - \theta} 
 \ \parallel \omega^\rho u \parallel_r^\theta \ .
 \eeq

We shall also use extensively the following Leibnitz estimates. \\

\noi {\bf Lemma 2.2.} {\it Let $1 < r,r_1,r_3 < \infty$ and
$$1/r = 1/r_1 + 1/r_2 = 1/r_3 + 1/r_4 \ .$$

\noi Then the following estimates hold for $\sigma \geq 0$~:}
\beq  
\label{2.2e}
\parallel \omega^\sigma (uv) \parallel_r \ \leq C \left ( \parallel  \omega^\sigma u \parallel_{r_1} 
 \ \parallel v \parallel_{r_2} + \parallel  \omega^\sigma v \parallel_{r_3} \  
 \parallel u \parallel_{r_4} \right )
\ .  \eeq 

An easy consequence of Lemmas 2.1 and 2.2 is the inequality 
\bea  
\label{2.3e}
\parallel \omega^\sigma f\ u \parallel \ &\leq& C  \left ( \parallel f \parallel_\infty \ + \ \parallel \omega^{n/2} f \parallel \right ) 
 \ \parallel  \omega^{\sigma} u \parallel  \  .\nn \\
&\leq& C \parallel \omega^{n/2\pm 0} f \parallel  \ \parallel  \omega^{\sigma} u \parallel 
 \eea

\noi which holds for $|\sigma | < n/2$. \par

Another consequence is the following lemma.\\

\noi {\bf Lemma 2.3.} {\it Let $0 < \sigma = \sigma_1  + \sigma_2$ and $\sigma_1 \vee \sigma_2 < n/2$. Then} 
\beq  
\label{2.4e}
\parallel \omega^{\sigma - n/2} (uv) \parallel \ \leq C  \parallel  \omega^{\sigma_1} u \parallel \ \parallel  \omega^{\sigma_2} v \parallel  \  .
 \eeq

We shall also need some commutator estimates, which are most conveniently stated in terms of homogeneous Besov spaces $\dot{B}_{r, q}^\sigma$ \cite{1r}. In the applications, we shall use only the fact that $\dot{B}_{2,2}^\sigma$ = $\dot{H}^\sigma$. The following lemma is an extension of Lemma~3.6 in \cite{5r} and may have independent interest. The proof will be given in Appendix A1.\\

\noi {\bf Lemma 2.4.} {\it Let $P_i$, $i = 1,2$ be homogeneous derivative polynomials of degree $\alpha_i$ or $\omega^{\alpha_i}$ for $\alpha_i \geq 0$. Let $\lambda > 0$. Then for any (sufficiently regular) functions $m$, $u$ and $v$ the following estimates hold.
\bea  
\label{2.5e}
&&|< P_1 u, [\omega^\lambda , m] P_2 v>| \leq  C \parallel  m ; \dot{B}_{r_0,2}^{\sigma_0} \cap \nabla^{-1} \omega^{1-\nu} L^{q_0} \parallel \ \parallel u; \dot{B}_{r_1,2}^{\sigma_1} \cap L^{q_1} \parallel \nn \\
&& \parallel v;  \dot{B}_{r_2,2}^{\sigma_2} \cap L^{q_2}\parallel 
 \eea
\noi with $0 \leq \nu \leq 1$, $1 \leq r_i, q_i \leq \infty$, $0 \leq i \leq 2$,
\beq  
\label{2.6e}
\delta (q_0) = \sigma_0 + \delta (r_0) - \nu \quad , \quad \delta (q_i) = \sigma_i + \delta (r_i) \ , \ i = 1,2 \ .
 \eeq
\beq  
\label{2.7e}
\sum_{0 \leq i \leq 2} \sigma_i + \delta (r_i) = \lambda + \alpha_1 + \alpha_2 + n/2
 \eeq
\beq  
\label{2.8e}
\left\{ \begin{array}{l} \sigma_0  + \left ( \sigma_1  \wedge \sigma_2 \right ) \geq \lambda + \alpha_1 + \alpha_2\\ \\ \sigma_1  +  \sigma_2  \geq \lambda + \alpha_1 + \alpha_2 - \nu \end{array} \right .
 \eeq
\noi where $\delta (r) \equiv n/2 - n/r$ and $\nabla^{-1} \omega^{1-\nu} L^q$ is the space of tempered distributions $m$ such that $\omega^{\nu - 1} \nabla m \in L^q$.}\\

\noi {\bf Remark 2.1.}  The condition (\ref{2.6e}) implies that the various spaces occuring in the RHS of (\ref{2.5e}) are homogeneous under dilation, and the condition (\ref{2.7e}) is the global homogeneity condition of the estimate. More general, possibly non homogeneous, estimates can be derived by the same method (see Appendix A1). \\

We shall repeatedly use the estimate of $s = \nabla \varphi$ with $\varphi$ defined by (\ref{1.18e})
\beq  
\label{2.9e}
\parallel \nabla^j s\parallel_\infty \ + \ \parallel  \omega^{n/2} \nabla^j s\parallel \ \leq \ C \parallel  \omega^{n/2\pm 0} \nabla^j s\parallel
\ \leq \ C \ t^{\lambda_j - 1} \parallel v_0; H^\rho \parallel^2 
 \eeq 

\noi for $j = 0,1$, where
\beq  
\label{2.10e}
\lambda_j = \gamma - (1/2) [1 + j + \gamma - 2 \rho ]_+ \ .
 \eeq

\noi The first inequality in (\ref{2.9e}) follows from Lemma 2.1 and the second one from 
$$\parallel \omega^{n/2\pm 0} \nabla^j s \parallel \ = \ C  \ t^{\gamma - 1}\parallel  \omega^{1+j+\gamma -n/2\pm 0} \chi_L |v_0|^2\parallel \ ,$$

\noi from the definition of $\chi_L$ and from Lemma 2.3. Note that up to an $\varepsilon$ in the case of equality
\beq  
\label{2.11e}
\lambda_0 = \gamma \wedge (1/2 + \delta ) \quad , \quad \lambda_1 = \gamma \wedge \delta
 \eeq

\noi where $\delta = \rho - 1 + \gamma /2$, so that $\lambda_1 > 0$ and $\lambda_0 > 1/2$ for $\gamma > 1/2$ and $\rho > 1 - \gamma /2$. The latter condition is the subcriticality condition mentioned in the introduction.\\

We shall also need some phase estimates. The following lemma is a variant of Lemma 3.3 in \cite{5r}. \\

\noi {\bf Lemma 2.5.} {\it Let $\varphi$ be a real function. Let $\sigma \geq 0$ and $1 \leq q, r\leq \infty$. Then the following estimate holds
$$\parallel\left  (\exp (i \varphi ) - 1\right ) ;  \dot{B}_{r,q}^\sigma\parallel \ \leq C  \ \parallel  \varphi ; \dot{B}_{r,q}^\sigma \parallel \left ( 1 + \ \parallel  \varphi ; \dot{B}_{\infty , \infty}^0\parallel  \right )^{[\sigma ]}$$

\noi where $[\sigma ]$ is the integral part of $\sigma$.} \\

An essential tool in this paper, as in \cite{4r} \cite{5r} \cite{6r}, is an estimate which follows from a local conservation law for solutions of a suitable linear Schr\"odinger equation.\\

\noi {\bf Lemma 2.6.} {\it Let $1/2 < \rho < n/2$, let $I$ be an interval, let $s \in L_{loc}^\infty (I, L^\infty \cap \dot{H}^{n/2})$, $s$ real ${I\hskip-1truemm R}^n$ vector valued, and let $v \in {\cal C}(I, H^\rho )$ be a solution of the equation
\beq
i \partial_t v + (1/2) \Delta_s v =  V v 
\label{2.12e}
\eeq

\noi in $I$ for some real $V \in L_{loc}^\infty(I, L^\infty )$. Then for any $t_1$, $t \in I$, $t_1 \leq t$, one can write 
\beq  
\label{2.13e}
|v(t)|^2 - |v(t_1)|^2 = V_1(t_1, t) + V_2(t_1, t)
 \eeq

\noi where $V_1$ and $V_2$ satisfy the following estimates~:
\beq  
\label{2.14e}
\parallel \omega^{2\sigma - 2 - n/2}  \ V_1(t_1, t) \parallel \ \leq \ C \int_{t_1}^t dt' \parallel  \omega^\sigma v(t') \parallel^2
\eeq 
 
\noi for $1/2 < \sigma \leq \rho \wedge (1 + n/4)$, and
\beq  
\label{2.15e}
\parallel \omega^{2\sigma - 1 - n/2}  \ V_2(t_1, t) \parallel \ \leq \ C \int_{t_1}^t dt' \left ( \parallel s(t') \parallel_\infty\ + \ \parallel  \omega^{n/2} s(t') \parallel \right ) \ \parallel  \omega^{\sigma} v(t') \parallel^2
\eeq 

\noi for $0 < \sigma \leq \rho$.}\\

\noi {\bf Sketch of proof.} The formal conservation law
\beq
\partial_t |v|^2 = - \ {\rm Im}\ \overline{v} \Delta_s v = - \ {\rm Im}\  \overline{v}  \Delta v + \nabla \cdot s |v|^2
\label{2.16e}
\eeq

\noi yields (\ref{2.13e}), where for any test function $\psi$ of the space variable 
\beq  
\label{2.17e}
<V_1(t_1, t), \psi > \ = - (i/2) \int_{t_1}^t dt' <v(t') , [\Delta , \psi ] v(t')>\ ,
 \eeq
\beq  
\label{2.18e}
<V_2(t_1, t), \psi > \ = - \int_{t_1}^t dt' <s(t') |v(t')|^2, \nabla \psi >\ .
 \eeq

\noi By Lemma 2.4 with $\alpha_i = 0$, $\lambda = 2$, $r_i = 2$ ($0 \leq i \leq 2$), $\sigma_1 = \sigma_2 = \sigma$ and therefore $\sigma_0 = n/2 + 2 - 2\sigma$, we obtain 
\beq  
\label{2.19e}
\left | <V_1(t_1, t), \psi > \right |\ \leq \ C \parallel \omega^{n/2 + 2 - 2\sigma}  \ \psi \parallel \int_{t_1}^t dt' \parallel \omega^\sigma  v(t') \parallel^2
\eeq

\noi for $1/2 < \sigma \leq \rho \wedge (1 + n/4)$, from which (\ref{2.14e}) follows by duality.\par

In order to estimate $V_2$, we estimate 
\beq
\label{2.20e}
\parallel \omega^{2\sigma - n/2}  \ s|v|^2 \parallel\ \leq \ C \left ( \parallel s\parallel_{\infty} \ + \ \parallel \omega^{n/2}  \ s \parallel \right ) \parallel \omega^\sigma  v \parallel^2
\eeq

\noi by (\ref{2.3e}) and Lemma 2.3. The estimate (\ref{2.15e}) then follows from (\ref{2.18e}) (\ref{2.20e}) by duality. \par\nobreak
\hfill $\sq$ \par

We next exploit the previous lemma in the relevant situation.\\

\noi {\bf Lemma 2.7.} {\it Let $1/2 < \rho < n/2$. Let $v_0 \in H^\rho$ and let $s = \nabla \varphi$ with $\varphi$ defined by (\ref{1.18e}). Let $I = (0, T]$ and let $v \in L^\infty (I, H^\rho ) \cap {\cal C} ([0, T], L^2)$ satisfy the equation (\ref{2.12e}) in $I$ for some real $V \in L_{loc}^\infty (I, L^\infty )$. Then $|v(t)|^2 - |v(0)|^2$ tends to zero in $\dot{H}^\mu$ for $-1 - n/2 < \mu < 2\rho - n/2$ when $t$ tends to zero. Furthermore
\beq
|v(t)|^2 - |v(0)|^2 = V_1(0,t) + V_2(0,t)
\label{2.21e}
\eeq

\noi with 
\beq  
\label{2.22e}
\parallel \omega^{2\sigma - 2 - n/2}  \ V_1(0, t) \parallel \ \leq \ C \ a^2\ t 
\eeq 
 
\noi for $1/2 < \sigma \leq \rho \wedge (1 + n/4)$, and
\beq  
\label{2.23e}
\parallel \omega^{2\sigma - 1 - n/2}  \ V_2(0, t) \parallel \ \leq \ C \ a^2\ a^2_0 \ t^{\lambda_0} 
\eeq 

\noi for $0 < \sigma \leq \rho$, with}
\beq  
\label{2.24e}
a \ = \ \parallel v; L^\infty (I, H^\rho )\parallel\quad , \quad a_0 \ = \ \parallel v_0 ; H^\rho \parallel \ .
 \eeq

\noi {\bf Proof.} We first prove that $V_1 (t_1, t)$ and $V_2 (t_1, t)$ defined by (\ref{2.17e}) (\ref{2.18e}) converge when $t_1 \to 0$ in the norms occuring in (\ref{2.14e}) (\ref{2.15e}) and that the limits satisfy (\ref{2.22e}) (\ref{2.23e}). This is obvious for $V_1$. As regards $V_2$, we estimate $s$ in (\ref{2.15e}) by (\ref{2.9e}) with $j = 0$. The resulting power of $t'$ in the integral is then integrable at $t= 0$ since $\lambda_0 > 0$. This proves the convergence of $V_2(t_1, t)$ as $t_1 \to 0$ and the estimate (\ref{2.23e}).\par

On the other hand, by Lemma 2.3, $|v(t)|^2$ is bounded in $\dot{H}^{2\rho - n/2}$ uniformly in $t$. Together with (\ref{2.13e}) and with the previous convergence of $V_1(t_1, t)$ and $V_2(t_1, t)$, this implies that $|v(t)|^2 - |v(t_1)|^2$ converges in $\dot{H}^\mu$ for $-1 - n/2 < \mu < 2\rho - n/2$ when $t_1$ tends to zero for fixed $t$. We next identify the limit. Now
$$|v(t)|^2 - |v(t_1)|^2 = |v(t)|^2 - |v(0)|^2 - \left ( |v(t_1)|^2 - |v(0)|^2 \right )$$

\noi and from Lemma 2.3
$$\parallel \omega^{\sigma - n/2}  \left (  |v(t_1) |^2  - |v(0)|^2 \right ) \parallel\ \leq \ C \parallel v(t_1) - v(0) \parallel \ \parallel \omega^{\sigma}  \left  ( v(t_1) + v(0) \right )\parallel$$

\noi for $0 < \sigma \leq \rho$, so that $|v(t)|^2 - |v(t_1)|^2$ tends to $|v(t)|^2 - |v(0)|^2$ in $\dot{H}^\mu$ for $-n/2 < \mu \leq \rho - n/2$ when $t_1$ tends to zero for fixed $t$. By an appropriate abstract argument, this implies that the same convergence holds in the whole range $-1 - n/2 < \mu < 2\rho - n/2$. This also implies (\ref{2.21e}), which together with the available estimates, completes the proof of the stated convergence. \par
\hfill $\sq$ \par

\noi {\bf Remark 2.2.} The difference $|v(t)|^2 - |v(0)|^2$ tends to zero in some norms which are not expected to be finite for $|v(t)|^2$ and $|v(0)|^2$ separately, typically in $\dot{H}^\mu$ for $-1-n/2 < \mu \leq - n/2$.\\

\noi {\bf Remark 2.3.} In most of the applications, we shall take $v_0 = v(0)$, but this is not needed in Lemma 2.7.

\mysection{The linearized Cauchy problem for v} 
\hspace*{\parindent}
In this section we study the Cauchy problem for the linearized equation (\ref{1.23e}) with $L(v)$ defined by (\ref{1.22e}) for a given $v$, with initial time $t_0 \geq 0$. We first give a preliminary result with $t_0 > 0$, where we do not study the behaviour of the solution as $t$ tends to zero. \\

\noi {\bf Proposition 3.1.} {\it Let $\rho > \gamma /2$, let $I = (0, T]$, let $v_0 \in H^\rho$ and let $v \in L_{loc}^\infty (I, H^\rho )$. Let $0 \leq \rho ' < n/2$, let $0 < t_0 \leq T$ and let $v'_0 \in H^{\rho '}$. Then the equation (\ref{1.23e}) has a unique solution $v' \in {\cal C} (I, H^{\rho '})$ with $v'(t_0) = v'_0$. The solution satisfies
$$\parallel v'(t) \parallel\ = \ \parallel v'_0 \parallel$$

\noi for all $t \in I$ and is unique in ${\cal C}(I, L^2)$.}\\

The proof is sketched in Appendix A2.\par

We next study the boundedness and continuity properties near $t=0$ of the solutions of (\ref{1.23e}) obtained in Proposition 3.1. Since we shall eventually be interested in taking $\rho ' = \rho$, we already impose the condition $\rho < n/2$ in the next proposition (see however Remark 3.2 below).\\

\noi {\bf Proposition 3.2.} {\it Let $1 -\gamma /2 < \rho < n/2$, let $I = (0, T]$ and let $v \in L^\infty (I, H^\rho ) \cap {\cal C} ([0, T], L^2)$ with $v(0) = v_0$. Let $s = \nabla \varphi$ with $\varphi$ defined by (\ref{1.18e}). Let $v$ satisfy the equation (\ref{2.12e}) in $I$ for some real $V \in L_{loc}^\infty (I, L^\infty )$. Let $1/2 \leq \rho ' < n/2$ and let $v' \in {\cal C} (I, H^{\rho '})$ be a solution of the equation (\ref{1.23e}) in $I$. Then \par

(1) $v' \in ({\cal C} \cap L^\infty ) (I, H^{\rho '}) \cap {\cal C}_w ([0, T] , H^{\rho '}) \cap {\cal C} ([0, T] , H^{\sigma})$ for $0 \leq \sigma < \rho '$.\par

(2) For all $t \in [0, T]$, $t_1 \in I$, the following estimate holds 
\beq  
\label{3.1e}
\parallel \omega^{\rho '}  v'(t) \parallel\ \leq \ \parallel \omega^{\rho '}  v'(t_1) \parallel\ E(|t- t_1|)
 \eeq

\noi where
\bea  
\label{3.2e}
E(t) &\equiv& E(t,a) = \exp \left \{ C \left ( a^2 t^{\lambda_1} + a^4  t^{2 \lambda_0 - 1} \right ) \right \} \ , \\
 a &=& \parallel v ; L^\infty (I, H^\rho ) \parallel
\label{3.3e}
 \eea
\noi and $\lambda_j$ is defined by (\ref{2.10e}).\par

(3) For all $t$, $t_1 \in [0, T]$, the following estimate holds} 
\beq  
\label{3.4e}
\parallel  v'(t) - v'(t_1)\parallel \ \leq C|t-t_1|^{(\rho '/2) \wedge (2\gamma - 1)} (1 + a^2)^2 \parallel  v'(t_1) ; H^{\rho '}\parallel \ .
 \eeq

\noi {\bf Remark 3.1.} The estimate (\ref{3.1e}) for $t$, $t_1 \in I$ holds for $0 \leq \rho ' < n/2$, as will be clear from the proof. The condition $\rho ' \geq 1/2$ is needed to derive (\ref{3.4e}) which is used in turn to extend (\ref{3.1e}) to $t = 0$.\\

\noi {\bf Remark 3.2.} The assumption $\rho < n/2$ in Proposition 3.2 can be dispensed with at the expense of using slightly different estimates, which yield different powers of $t$ in (\ref{3.2e})  and (\ref{3.4e}).\\

\noi {\bf Proof.} We know already that the $L^2$- norm of $v'$ is conserved. The bulk of the proof consists in deriving the estimates (\ref{3.1e}) and (\ref{3.4e}) for $t$, $t_1 \in I$. We begin with (\ref{3.1e}). From (\ref{1.22e}) (\ref{1.23e}) we obtain 
\bea  
\label{3.5e}
\partial_t \parallel \omega^{\rho '} v' \parallel^2 &=&  {\rm Im} \ <v', [\omega^{2\rho '}, L(v)]v'>\nn \\
&=& \ {\rm Re} \ <v', [\omega^{2\rho '}, s] \cdot \nabla v'>\ + \ {\rm Im} \ <v', [\omega^{2\rho '}, f] v'>
 \eea

\noi where 
\beq  
\label{3.6e}
f = (1/2) s^2 + t^{\gamma - 2} g_S (v) + \left ( t^{\gamma - 2} g_L (v_0) - \partial_t \varphi \right ) + t^{\gamma - 2} \left ( g_L (v) - g_L(v_0) \right ) \ .
 \eeq

\noi We estimate the first term in the RHS of (\ref{3.5e}) by Lemma 2.4 with $\lambda = 2\rho '$, $\alpha_1 = 0$, $\alpha_2 = 1$, $r_i = 2$, $\sigma_1 = \sigma_2 = \rho '$, so that $\sigma_0 = 1 + n/2$ and $q_0 = \infty$. \par

We estimate similarly the last term by Lemma 2.4 with $\lambda =2 \rho '$, $\alpha_1 = \alpha_2 = 0$, $r_i = 2$, $\sigma_1 = \sigma_2 = \rho '$, so that $\sigma_0 = n/2$ and $\delta (q_0) = n/2-1$. We obtain
\beq  
\label{3.7e}
\left |\partial_t \parallel \omega^{\rho '} v' \parallel^2 \right | \ \leq \ C\left ( \parallel \omega^{n/2} \nabla s \parallel\ + \ \parallel \nabla s \parallel_\infty \ + \ \parallel \omega^{n/2} f \parallel \right ) \parallel \omega^{\rho '} v' \parallel^2 \ .
\eeq

\noi We estimate the various norms successively. We first estimate $\nabla s$ by (\ref{2.9e}) with $j = 1$ so that
\beq  
\label{3.8e}
\parallel \omega^{n/2} \nabla s \parallel\ + \ \parallel \nabla s \parallel_\infty \ \leq \ C\ a^2\ t^{\lambda_1 - 1}
 \eeq

\noi and similarly
\beq  
\label{3.9e}
\parallel \omega^{n/2} s^2 \parallel\  \leq \ \parallel \omega^{n/2\pm 0} s \parallel^2\ \leq \ C\ a^4 \ t^{2\lambda_0 - 2}
 \eeq

\noi by Lemma 2.1 and (\ref{2.9e}) with $j = 0$.\par

We next estimate 
\bea  
\label{3.10e}
t^{\gamma - 2} \parallel \omega^{n/2} g_S(v) \parallel &\leq& C\ t^{\gamma - 2+ \rho - \gamma /2} \parallel \omega^\rho v \parallel^2\nn \\
&\leq &C\ a^2\ t^{\lambda_1 - 1} 
 \eea	

\noi for $\rho \geq \gamma /2$, 	and similarly (see (\ref{1.19e}))
\bea  
\label{3.11e}
&&\parallel \omega^{n/2} \left ( \partial_t \varphi - t^{\gamma - 2}  g_L(v_0)\right ) \parallel \ = \ (1 - \gamma )^{-1} \ t^{\gamma - 2} \parallel \omega^{n/2} \ \widetilde{\chi}_L g(v_0)\parallel \nn \\
&&\leq C\ t^{\gamma - 2+ \rho - \gamma /2} \parallel \omega^\rho v_0 \parallel^2 \ \leq \ C\ a^2\ t^{\lambda_1 - 1} \ .
 \eea	

\noi The contribution of the last term in $f$ is estimated by the use of Lemma 2.7. From  (\ref{2.21e}), from  (\ref{2.22e}) with $\sigma = (1 + \gamma /2) \wedge \rho > 1/2$ and from (\ref{2.23e}) with\break\noindent $\sigma = (1 + \gamma )/2 \wedge \rho > 0$, we obtain
\bea  
\label{3.12e}
t^{\gamma - 2} \parallel \omega^{n/2} \left ( g_L(v) - g_L(v_0 \right )) \parallel &=& C\   t^{\gamma - 2} \parallel \omega^{\gamma - n/2}\chi_L \left ( V_1(0,t) + V_2 (0, t)\right )  \parallel\nn \\
&\leq &C\left ( a^2\ t^{\lambda_1 - 1} + a^4\ t^{2\lambda_0 - 2} \right ) \ .
 \eea	

\noi Collecting  (\ref{3.7e})- (\ref{3.12e}), we obtain
\beq  
\label{3.13e}
\left | \partial_t \parallel \omega^{\rho '} v' (t) \parallel^2\right |  \ \leq \ N(t) \parallel \omega^{\rho '} v'(t) \parallel^2 
\eeq

\noi where
\beq  
\label{3.14e}
N(t) = C\left ( a^2\ t^{\lambda_1 - 1} + a^4 \ t^{2\lambda_0 - 2}\right ) \ .
 \eeq	

\noi The crucial point of this estimate is that $N(t)$ is integrable in time at $t=0$ since $\gamma > 1/2$ and $\rho > 1 - \gamma /2$ (see  (\ref{2.11e})). Note that 
$$\gamma /2 < 1/2 < 1 - \gamma /2 < (1 + \gamma )/2 < 1 + \gamma /2$$

\noi for $1/2 < \gamma < 1$. Therefore the condition $\rho > 1 - \gamma /2$ implies the conditions $\rho \geq \gamma /2$ and $\rho > 1/2$ used in the proof of  (\ref{3.10e}) and  (\ref{3.12e}) respectively. Furthermore, there exists an interval, namely $1 - \gamma /2 < \rho < (1 + \gamma )/2$ where the [ ]$_+$ brackets are inactive. The estimate  (\ref{3.1e})  (\ref{3.2e}) follows from  (\ref{3.13e})  (\ref{3.14e}) by integration for $t_1$, $t\in I$. \par

We next derive the estimate (\ref{3.4e}) for $t$, $t_1 \in I$. For that purpose we define (see  (\ref{1.8e}))
\beq  
\label{3.15e}
\widetilde{v}'(t) = U(-t) v'(t)
 \eeq	
\beq  
\label{3.16e}
\widetilde{L} = L(v) + (1/2) \Delta = is \cdot \nabla + (i/2) (\nabla \cdot s ) + f
 \eeq

\noi with $f$ given by (\ref{3.6e}). We rewrite (\ref{1.23e}) as
\beq  
\label{3.17e}
i \partial_t \widetilde{v}'  = U(-t) \widetilde{L} U(t) \widetilde{v}'
 \eeq	

\noi so that for $t$, $t_1 \in I$, for fixed $t_1$,
\bea
\label{3.18e}
\partial_t \parallel \widetilde{v}' (t) - \widetilde{v}' (t_1) \parallel^2 &=&2 \ {\rm Im} \ < \widetilde{v}' (t) - \widetilde{v}' (t_1) , U(-t) \widetilde{L} \ U(t) \  \widetilde{v}' (t_1) >\nn \\
&=&2 \ {\rm Im} \ < w, \widetilde{L}\  v_*>
\eea

\noi where
\beq
\label{3.19e}
\left \{ \begin{array}{l} v_* = U(t-t_1) v'(t_1) \\ \\ w = v'(t) - v_* \ . \end{array} \right . 
\eeq

\noi We estimate 
\bea  
\label{3.20e}
&&\left |\partial_t \parallel w \parallel^2 \right | \ \leq \ 2\left | {\rm Re} \ <w, s\cdot \nabla v_*> \right |\nn \\
&&+\ C  \parallel w\parallel \left ( \parallel \omega^{n/2- \rho_1} \nabla \cdot s \parallel\ + \ \parallel \omega^{n/2- \rho_1} f \parallel \right )  \parallel \omega^{\rho_1} v'(t_1) \parallel 
\eea

\noi for some $\rho_1$ with $0 < \rho_1 \leq \rho '$, to be chosen later. \par

For $0 < \rho ' < 1$, we write
\beq
\label{3.21e}
<w,s\cdot \nabla v_*> \ = \ - <\omega^{- \rho '} \nabla \cdot sw, \omega^{\rho '} v_*>
\eeq

\noi and we estimate by Lemma 2.2
\beq
\label{3.22e}
\left |< w, s\cdot \nabla v_*>\right |  \ \leq \ C \parallel \omega^{1 - \rho '} w \parallel\left (  \parallel s \parallel_\infty \ + \ \parallel \omega^{n/2} s \parallel \right ) \ \parallel \omega^{\rho '} v'(t_1) \parallel \ .
\eeq

\noi For $\rho ' = 1$, we estimate
\beq
\label{3.23e}
\left |< w, s\cdot \nabla v_*>\right |  \ \leq \  \parallel w \parallel\ \parallel s \parallel_\infty \  \parallel \omega^{\rho '} v'(t_1) \parallel  \ .
\eeq

\noi For $\rho ' > 1$, we estimate
\beq
\label{3.24e}
\left |< w, s\cdot \nabla v_*>\right |  \ \leq \ C \parallel w \parallel\ \parallel \omega^{n/2- \rho_1} \nabla s \parallel  \ \parallel \omega^{\rho _1} v'(t_1) \parallel 
\eeq

\noi for $1 < \rho_1 \leq \rho '$.\par

Collecting (\ref{3.20e})-(\ref{3.24e}) yields 
\bea  
\label{3.25e}
&&\left |\partial_t \parallel w \parallel^2 \right | \ \leq \ C\Big \{ \chi  (\rho ' \leq 1) \parallel \omega^{1 - \rho '} w \parallel\left (  \parallel s \parallel_\infty \ + \ \parallel \omega^{n/2} s \parallel\right ) \parallel \omega^{\rho '} v'(t_1) \parallel \nn \\
&&+\   \parallel w\parallel \left ( \parallel \omega^{n/2- \rho_1} \nabla s \parallel\ + \ \parallel \omega^{n/2- \rho_1} f \parallel \right )  \parallel \omega^{\rho_1} v'(t_1) \parallel \Big \}
\eea

\noi with $0 < \rho_1 \leq \rho '$ and $\rho_1 > 1$ in the $\nabla s$ term if $\rho ' > 1$.\par

For $1/2 \leq \rho ' \leq 1$, we interpolate 
$$\parallel \omega^{1 - \rho '} w \parallel \ \leq \ y^\theta \parallel \omega^{\rho '} w \parallel^{1/\rho ' - 1}$$

\noi where
$$y = \ \parallel w(t) \parallel^2 \qquad , \quad \theta = 1 - 1/(2 \rho ')$$

\noi so that (\ref{3.25e}) becomes
\bea  
\label{3.26e}
\left |\partial_t y \right | &\leq& \ C\left \{ \chi  (\rho ' \leq 1) \left ( \parallel s \parallel_\infty \ + \ \parallel \omega^{n/2} s \parallel\right ) a{'}_1^{1/\rho '} y^\theta \right . \nn \\
&&\left . +\   \left ( \parallel \omega^{n/2- \rho_1} \nabla s \parallel \ + \ \parallel \omega^{n/2- \rho_1} f \parallel \right )  a'_1 \ y^{1/2}  \right \} 
\eea

\noi with $a'_1 = \parallel v'(t_1);H^{\rho '}\parallel$. We estimate $s$ in the first term in the RHS of (\ref{3.26e}) by (\ref{2.9e}) with $j = 0$ and we estimate the various contributions to the second term for suitable values of $\rho_1$. We first estimate
\bea
\label{3.27e}
\parallel \omega^{n/2- \rho_1} \nabla s \parallel &\leq& C\ t^{\gamma - 1} \parallel \omega^{2 + \gamma - \rho_1 - n/2} \chi_L |v_0|^2 \parallel \nn \\
&\leq& C\ t^{\gamma - 1} \parallel \omega^{\rho_2} v_0  \parallel^2
\eea

\noi by Lemma 2.3 with $0 < \rho_1 \leq \rho '$, $0 < \rho_2 \leq \rho$, $\rho_1 > 1$ if $\rho ' > 1$ and $\rho_1 + 2 \rho_2 = 2 + \gamma$, in the case where $\rho ' + 2 \rho \geq 2 + \gamma$. In the opposite case, we take $\rho_1 = \rho '$, $\rho_2 = \rho$, and we use the cut off $\chi_L$ so that finally, for the relevant choice of $\rho_1$, 
\beq
\label{3.28e}
\parallel \omega^{n/2 - \rho_1} \nabla s \parallel \ \leq \ C\ a^2 \ t^{\mu_1 - 1} 
\eeq

\noi with
\beq
\label{3.29e}
\mu_1 = \gamma - (1/2) (2 + \gamma - 2\rho - \rho ')_+  \ .
\eeq

\noi Similarly, we estimate
\bea
\label{3.30e}
&&\parallel \omega^{n/2- \rho_1} |s|^2 \parallel \ \leq \ C\parallel \omega^{n/2 - \rho_1/2} s \parallel^2 \nn \\
&&= \  C\ t^{2\gamma - 2} \parallel \omega^{1 + \gamma - \rho_1/2 - n/2} \chi_L |v_0|^2 \parallel^2 \ \leq \ C\ a^4\ t^{2\mu_0 - 2}
\eea

\noi with
\beq
\label{3.31e}
\mu_0 = \gamma - (1/2) (1 + \gamma - 2\rho - \rho '/2)_+  \ .
\eeq

\noi We next estimate, with $\rho_1 = \rho '$
\bea
\label{3.32e}
&&t^{\gamma -2}  \parallel \omega^{n/2 - \rho '} g_S(v) \parallel \ \ \leq \ C\ a^2\  t^{\mu_1 -1}\ , \\
&&\parallel \omega^{n/2- \rho '} \left ( \partial_t \varphi - t^{\gamma -2} g_L(v_0) \right ) \parallel \ \leq \ C\ a^2 \ t^{\mu_1 - 1} \ . 
\label{3.33e}
\eea

\noi We next consider
$$t^{\gamma -2}  \parallel \omega^{n/2 - \rho_1} \left ( g_L(v) - g_L(v_0)\right )  \parallel \ = \  C\  t^{\gamma -2} \parallel \omega^{\gamma - \rho_1 - n/2} \chi_L \left ( V_1 + V_2 \right ) \parallel$$

\noi where $V_i \equiv V_i (0, t)$, $i = 1,2$, are defined by (\ref{2.17e}) (\ref{2.18e}). By (\ref{2.14e}) we estimate
\beq
\label{3.34e}
\parallel \omega^{\gamma - \rho_1 - n/2 } \chi_L \ V_1\parallel \ \leq \  C \int_0^t dt' \parallel \omega^{\rho_2}  v(t') \parallel^2
\eeq

\noi for $0 < \rho_1 \leq \rho '$, $1/2 < \rho_2 \leq \rho$ and $\rho_1 + 2\rho_2 = 2 + \gamma$ in the case where $\rho ' + 2 \rho \geq 2 + \gamma$. In the opposite case, we take $\rho_1 = \rho '$, $\rho_2 = \rho$ and we use the cut off $\chi_L$, so that finally, for the relevant choice of $\rho_1$,
\beq
\label{3.35e}
t^{\gamma -2}  \parallel \omega^{\gamma - \rho_1 - n/2} \chi_L  \ V_1  \parallel \ \leq \  C\  a^2\ t^{\mu_1 -1} \ .
\eeq

\noi We next estimate
\bea
\label{3.36e}
&&\parallel \omega^{\gamma - \rho_1 - n/2 } \chi_L \ V_2\parallel \ \leq \  \int_0^t dt' \parallel \omega^{1 + \gamma - \rho_1- n/2}  \ \chi_L(t) \left  (s|v|^2\right ) (t') \parallel \nn \\
&&\leq \ C \int_0^t dt' \parallel \omega^{1 + \gamma - \rho_1/2- n/2} |v(t')|^2\parallel \ \parallel  \omega^{n/2 - \rho_1/2} s(t') \parallel \nn \\
 &&\leq \ C  \int_0^t dt' \ t{'}^{\gamma -1} \parallel \omega^{\rho_2}   v(t')  \parallel^2 \  \parallel \omega^{1 + \gamma - \rho_1/2- n/2} \ \chi_L(t')   |v_0|^2 \parallel  \nn \\
&&\leq\  C  \parallel \omega^{\rho_2}  v _0\parallel^2 \int_0^t dt'  t{'}^{\gamma -1} \parallel \omega^{\rho_2}  v(t')  \parallel^2 
\eea

\noi by repeated use of Lemma 2.3, for $0 < \rho_1 \leq \rho '$, $\rho_1 < 1 + \gamma$, $0 < \rho_2 \leq \rho$ and $\rho_1/2 + 2 \rho_2 = 1 + \gamma$, in the case where $\rho ' /2 + 2 \rho \geq 1 + \gamma$. In the opposite case, we take $\rho_1 = \rho '$, $\rho_2 = \rho$ and we use the cut offs $\chi_L(t)$ and $\chi_L(t')$, so that finally, for the relevant choice of $\rho_1$,  
\bea
\label{3.37e}
&&t^{\gamma -2}  \parallel \omega^{\gamma - \rho_1 - n/2} \chi_L\  V_2\parallel \ \leq \ C\parallel v_0; H^\rho\parallel^2 t^{\mu_0  -2} \nn \\
&&\times \int_0^t dt' \ t{'}^{\mu_0  -1} \parallel v(t') ; H^\rho \parallel^2 \ \leq \  C\  a^4\ t^{2\mu_0 -2} \ .
\eea

\noi (Note that in the second case
$$\rho ' < 2 + 2\gamma - 4 \rho < 4 \gamma - 2 \leq 1 + \gamma \qquad \hbox{for $\gamma \leq 1$)}\ .$$

\noi Collecting (\ref{3.26e}), (\ref{2.9e}) with $j = 0$ and (\ref{3.28e}) (\ref{3.30e}) (\ref{3.32e}) (\ref{3.33e}) (\ref{3.35e}) (\ref{3.37e}) yields
\bea  
\label{3.38e}
\left |\partial_t y \right | &\leq&  C\left \{ \chi  (\rho ' \leq 1)  a^2 \ t^{\lambda_0 - 1} \ a{'}_1^{1/\rho '} \ y^\theta \right . \nn \\
&&\left . +   \left (  a^2\ t^{\mu_1 -1} + a^4 \   t^{2\mu_0 -2}\right )   a'_1\ y^{1/2} \right \} \ .
\eea

\noi Using the fact that the differential inequality
$$\left |\partial_t y \right |  \leq \ \sum_i b_i \ t^{\nu_i - 1}\ y^{\theta_i}$$

\noi with $0 \leq \theta_i < 1$, $\nu_i > 0$ implies
$$y(t) \leq \ C \sum_i \left ( b_i \ \nu_i^{-1} \left | t^{\nu_i} - t_1^{\nu_i}\right | \right )^{1/(1 - \theta_i)}$$

\noi for $t$, $t_1 > 0$ and $y(t_1) = 0$, we obtain
\bea  
\label{3.39e}
\parallel w \parallel &\leq&  C\left \{ \chi  (\rho ' \leq 1)  a^{2\rho '}   | t - t_1|^{\rho ' \lambda_0} + a^2   | t - t_1|^{\mu_1} + a^4 |t-t_1|^{2\mu_0-1}  \right \}  a{'}_1\nn \\
&\leq& C \left ( \chi (\rho' \leq 1) a^{2\rho '}+ a^2 (1 + a^2)\right ) |t-t_1| ^\mu \ a'_1  \ .
\eea

\noi where
\bea  
\label{3.40e}
\mu &=& \rho ' \lambda_0 \wedge \mu_1 \wedge (2 \mu_0 - 1) \nn \\
&=& \rho ' \gamma \wedge \rho ' (1/2 + \delta ) \wedge \gamma \wedge (\rho ' /2 + \delta ) \wedge (2\gamma - 1) \wedge (\rho ' /2 + 2\delta )\nn \\
&\geq& \rho ' /2 \wedge (2 \gamma - 1)
\eea

\noi since $1/2 < \gamma < 1$ and $\delta \equiv \rho - 1 + \gamma /2 > 0$. \par

On the other hand, we estimate
\bea  
\label{3.41e}
&&\parallel v'(t) - v'(t_1) \parallel \ \leq \ \parallel w \parallel \ + \ \parallel \left ( U(t-t_1) - 1\right ) v'(t_1) \parallel \nn \\
&&\leq \ \parallel w \parallel\ + \ |t-t_1|^{(\rho '/2) \wedge 1} \parallel \omega^{\rho ' \wedge 2} \ v'(t_1) \parallel \ .
\eea

\noi Collecting (\ref{3.39e}) (\ref{3.41e}) yields (\ref{3.4e}) for $t$, $t_1 \in I$. \par

We now exploit (\ref{3.1e}) and (\ref{3.4e}) in $I$ to complete the proof of the proposition. From (\ref{3.1e}) it follows that $v' \in L^\infty (I, H^{\rho '})$. From (\ref{3.1e}) and  (\ref{3.4e}) it then follows that $v'$ has a limit $v'(0)$ in $L^2$ and that (\ref{3.4e}) holds for $t$, $t_1 \in [0, T]$. It then follows by a standard abstract argument that $v'(0) \in H^{\rho '}$, that $v' \in {\cal C}_w ([0, T], H^{\rho '}) \cap {\cal C} ([0, T], H^{\sigma})$ for $0 \leq \sigma < \rho '$, and that (\ref{3.1e}) holds for all $t \in [0, T]$, $t_1 \in I$. \par\nobreak \hfill $\sq$\par

We have not proved so far that $v' \in {\cal C} ([0, T], H^{\rho '})$. This is true but requires a separate argument. \\

\noi {\bf Proposition 3.3.} {\it Under the assumptions of Proposition 3.2, $v' \in {\cal C} ([0, T], H^{\rho '})$ and (\ref{3.1e}) holds for all $t$, $t_1 \in [0, T]$.}\\

\noi The proof is identical with that of Proposition 3.3 of \cite{4r}. \par

We can now state the main result on the Cauchy problem for the linearized equation (\ref{1.23e}). \\

\noi {\bf Proposition 3.4.} {\it Let $1 - \gamma /2 < \rho < n/2$. Let $I = (0, T]$ and let $v \in L^\infty (I, H^\rho ) \cap  {\cal C} ([0, T], L^2)$ with $v(0) = v_0$. Let $s = \nabla \varphi$ with $\varphi$ defined by (\ref{1.18e}). Let $v$ satisfy the equation (\ref{2.12e}) in $I$ for some real $V \in L_{loc}^\infty (I, L^\infty )$. Let $1/2 < \rho ' < n/2$ and let $v'_0 \in H^{\rho '}$. Let $t_0 \in [0, T]$. Then there exists a unique solution $v' \in  {\cal C} ([0, T], H^{\rho '})$ of the equation (\ref{1.23e}) with $v'(0) = v'_0$. Furthermore $v'$ satisfies the estimates (\ref{3.1e}) and (\ref{3.4e}) for all $t$, $t_1 \in [0, T]$. The solution is actually unique in ${\cal C} ([0, T], L^2)$.}\\

\noi {\bf Proof.} For $t_0 > 0$, the result follows from Propositions 3.1, 3.2 and 3.3. For $t_0 =0$, it will be proved by a limiting procedure on $t_0$. For any $t_1 \in I$, let $v'_{t_1}$ be the solution of (\ref{1.23e}) with $v'_{t_1}(t_1) = v'_0$ given by Propositions 3.1 and 3.2. Let now $0 < t_1 <t_2 \leq T$. It follows from (\ref{3.1e}) that 
\beq
\label{3.42e}
\parallel \omega^{\rho '} v'_{t_i} (t)\parallel \ \leq \ E \left ( |t - t_i| \right ) \parallel \omega^{\rho '} v'_0\parallel 
\eeq

\noi for $i = 1,2$ and for all $t \in [0, T]$. Furthermore, from  (\ref{3.4e}) and  (\ref{3.42e}) and from $L^2$-norm conservation, it follows that
\begin{eqnarray}
\label{3.43e}
&&\parallel v'_{t_2}(t) - v'_{t_1}(t) \parallel  \ = \  \parallel  v'_{t_2}(t_1) - v'_0\parallel \ = \  \parallel  v'_{t_2}(t_1)  - v'_{t_2}(t_2) \parallel\nn  \\
&&\leq \ C|t_2 - t_1|^{(\rho '/2)\wedge (2\gamma - 1)} (1 + a^2)^2 \parallel v'_0 ; H^{\rho '}\parallel \ .
\end{eqnarray}

\noi From  (\ref{3.43e}) it follows that $v'_{t_1}$ converges in $L^\infty (I, L^2)$-norm to some $v' \in {\cal C}([0,T], L^2)$ when $t_1 \to 0$. From the uniform estimate (\ref{3.42e}) it follows by abstract arguments that $v' \in ({\cal C}_w \cap L^\infty )([0,T], H^{\rho '}) \cap {\cal C} ([0, T], H^\sigma )$ for $0 \leq \sigma < \rho '$, that $v'$ satisfies the estimates of Proposition 3.2 and that $v'(0) = v'_0$. Furthermore $v'$ is easily seen to satisfy (\ref{1.23e}) in $I$, so that $v' \in {\cal C}(I, H^{\rho '})$. It remains to be proved that actually $v'$ is strongly continuous in $H^{\rho '}$ at $t = 0$. This follows from Proposition 3.3, which has not been used so far. Alternatively it follows from the estimate (\ref{3.42e}) with $t_1 = 0$ that 
$$\lim_{t\to 0} \ \sup \parallel \omega^{\rho '}v'(t)\parallel\ \leq \ \parallel \omega^{\rho '}v'_0\parallel\ E(0) = \ \parallel \omega^{\rho '}v'(0)\parallel$$

\noi which together with weak continuity implies strong continuity at $t=0$. \par\nobreak \hfill $\sq$

\noi {\bf Remark 3.3.} Note that in the case where $t_0 = 0$, Proposition 3.3 is not needed for the proof of Proposition 3.4.

\mysection{The nonlinear Cauchy problem at time zero for v and u$_{\bf c}$}
\hspace*{\parindent} 
In this section we prove that the nonlinear equation (\ref{1.21e}) for $v$ with initial data at time $t_0$ has  a unique solution in a small time interval. We then rewrite that result in terms of $u_c$, related to $v$ by (\ref{1.12e}), and we give some additional bounds and regularity properties for $u_c$. In order to solve the equation (\ref{1.21e}) for $v$, we show that the map $\Gamma : v \to v'$ defined by Proposition 3.4 with $t_0 = 0$ is a contraction. For that purpose, we need to estimate the difference of two solutions of the linearized equation (\ref{1.23e}). For any pair of functions or operators $(f_1, f_2)$, we define
$$f_\pm = (1/2) \left ( f_2 \pm f_1 \right ) \ .$$
\vskip 3 truemm

\noi {\bf Lemma 4.1.} {\it Let $1 - \gamma /2 < \rho < n/2$. Let $I = (0,T]$ and let $v_i$, $i = 1,2$ satisfy the assumptions of Proposition 3.4 with $v_i(0) = v_0$. Let $1/2 < \rho ' < n/2$ and let $v'_i$, $i =1,2$ be the solutions of the equation (\ref{1.23e}) with $v'_i (0) = v'_0 \in H^{\rho '}$ obtained in Proposition 3.4. Then the following estimate holds for all $t$, $0 < t \leq T$~:
 \beq
\label{4.1e}
\parallel v'_-; L^\infty ((0, t], H^{\rho '}) \parallel \ \leq \ C\ E(t,a) aa' \left ( t^{\lambda_1}  + a^2 \ t^{2\lambda_0 - 1}\right )  \parallel v_-; L^\infty ((0, t], H^\rho ) \parallel
\eeq 

\noi where $E(t, a)$ is defined by (\ref{3.2e}) and}
 \beq
\label{4.2e} 
a = \ {\rm Max}\parallel v_i; L^\infty (I, H^\rho ) \parallel\ , \quad a' = \ {\rm Max}\parallel v'_i; L^\infty (I, H^{\rho '} ) \parallel\  .
\eeq

\noi {\bf Proof.} From (\ref{1.23e}) we obtain
$$i\partial_t v'_- = L_2 \ v'_- + L_- \ v'_1$$

\noi where $L_i =L(v_i)$, $g_i = g(v_i)$, so that
$$L_- = t^{\gamma - 2}\ g_-  \ .$$

\noi We estimate for $0 \leq \sigma  \leq \rho '$
\beq
\label{4.3e}
\partial_t \parallel \omega^{\sigma } v'_- \parallel^2 \ = \ 2\ {\rm Im}  \left (  < \omega^{\sigma }  v'_-, \omega^{\sigma } L_2 v'_- > \ + \ <   \omega^{\sigma } v'_-, \omega^{\sigma } L_- v'_1 > \right )\ .
\eeq 

\noi By the estimates in the proof of Proposition 3.2 (see in particular  (\ref{3.1e}) ; see also Remark 3.1), we obtain
\beq
\label{4.4e}
\parallel \omega^{\sigma } v'_- (t)\parallel \ \leq \ E(t,a) \int_0^t dt'\ t{'}^{\gamma -2} \parallel  \omega^{\sigma } g_- \ v'_1(t')\parallel \ .
\eeq

\noi We next estimate
 \bea
\label{4.5e}
&&\parallel \omega^{\sigma } g_- \ v'_1 \parallel \ \leq \ C \parallel \omega^{n/2 \pm 0} g_- \parallel \ \parallel \omega^\sigma  v'_1\parallel  \ , \\
&&\parallel \omega^{n/2 \pm 0} g_{S-}  \parallel \ \leq \ C \ t^{\rho - \gamma /2} \parallel  \omega^{\rho} v_-\parallel \ \parallel  \omega^{\rho} v_+\parallel\ .
\label{4.6e}
\eea 

\noi In order to estimate $g_{L-}$, we use again Lemmas 2.6 and 2.7. From the conservation law (\ref{2.16e}) and from the fact that $v_- (0) = 0$ we obtain (see (\ref{2.21e}))
\beq
\label{4.7e}
\left ( |v(t)|^2 \right )_- = V_{1-}(0,t) + V_{2-}(0,t)
\eeq

\noi where (see (\ref{2.17e}) (\ref{2.18e}))
\begin{eqnarray*}
&&V_{1-}(0,t) = - \int_0^t dt'\ {\rm Im} \left ( \overline{v}_+ \Delta v_- + \overline{v}_- \Delta v_+\right ) (t') \ , \\
&& V_{2-}(0,t) = \int_0^t dt'\ \nabla \cdot \left ( 2s\  {\rm Re}\ \overline{v}_+ v_- \right ) (t') \ .
\end{eqnarray*}

\noi By the same estimates as in Lemma 2.6, we obtain
\beq
\label{4.8e}
\parallel \omega^{2\sigma - 2 -n/2} V_{1-}(0,t) \parallel \ \leq \ C \int_0^t dt' \left ( \parallel \omega^{\sigma} v_{+} \parallel \ \parallel \omega^{\sigma} v_{-} \parallel\right ) (t') \eeq

\noi for $1/2 < \sigma \leq \rho \wedge (1 + n/4)$,
\beq
\parallel \omega^{2\sigma - 1 -n/2} V_{2-}(0,t) \parallel \ \leq \ C \int_0^t dt' \left ( \left ( \parallel s\parallel_\infty \ + \ \parallel \omega^{n/2} s \parallel \right ) \parallel \omega^{\sigma} v_{+} \parallel \ \parallel \omega^{\sigma} v_{-} \parallel\right ) (t')
\label{4.9e}
\eeq

\noi for $0 < \sigma \leq \rho$. In the same way as in Proposition 3.2 (see especially (\ref{3.12e})), we obtain 
\beq
\label{4.10e}
t{'}^{\gamma - 2} \parallel \omega^{\sigma } g_- \ v'_1 \parallel \ \leq \ C \ aa' \left ( t{'}^{\lambda_1 - 1} + a^2 \ t{'}^{2\lambda_0 - 2 }\right ) \parallel v_-; L^\infty ((0, t'], H^{\rho} ) \parallel \ . 
\eeq

\noi Substituting (\ref{4.10e}) into (\ref{4.4e}) yields (\ref{4.1e}). \par\nobreak \hfill $\sq$ \par

We can now state the main result on the Cauchy problem at time zero for the equation (\ref{1.21e}). \\

\noi {\bf Proposition 4.1.} {\it Let $1 - \gamma /2 < \rho < n/2$, let $v_0 \in H^\rho$ and define $\varphi$ by (\ref{1.18e}). Then there exists $T > 0$ and there exists a unique solution $v \in {\cal C}([0, T], H^\rho )$ of the equation (\ref{1.21e}) with $v(0) = v_0$. One can ensure that 
  \bea
\label{4.11e}
&& \parallel v; L^\infty ([0, T], H^\rho )  \parallel \ \leq \ R=2 \parallel v_0; H^\rho \parallel \\
&&C\ R^2 \left ( T^{\lambda_1} + R^2 T^{2\lambda_0 - 1} \right ) = 1
\label{4.12e}
\eea

\noi for some $C$ independent of $v_0$.}\\

\noi {\bf Proof.} Let $T > 0$. Let $F(T, v_0)$ be the set of $v \in  {\cal C}([0, T], H^\rho )$ such that $v(0) = v_0$ and satisfying the equation (\ref{2.12e}) in $(0, T]$ for some real $V \in L_{loc}^\infty ((0, T], L^\infty )$. It follows from Proposition 3.4 that $F(T, v_0)$ is stable under the map $\Gamma : v\to v'$ defined by that proposition with $t_0 = 0$ and $v'_0 = v_0$. In fact, $v'$ satisfies the equation (\ref{2.12e}) with
$$V = t^{\gamma -2} g(v) - \partial_t \varphi \in L_{loc}^\infty ((0, T], L^\infty )\ .$$

\noi Let $B(R)$ be the ball of radius $R$ in ${\cal C}([0, T], H^\rho )$. From Proposition 3.2 it follows that $B(R) \cap F(T, v_0)$ is stable under $\Gamma$ if
\beq
\label{4.13e}
E(T, R) \leq 2
\eeq 

\noi with $R = \ 2 \parallel v_0; H^\rho \parallel$. Furthermore by Lemma 4.1, $\Gamma$ is a contraction in the $L^\infty ([0, T]), H^\rho )$-norm on that set under the condition (\ref{4.12e}) for a suitable $C$. Such a condition at the same time implies (\ref{4.13e}). Therefore for $T$ sufficiently small to satisfy (\ref{4.12e}), the map $\Gamma$ has a unique fixed point in $B(R)$ provided $F(T,v_0)$ is non empty. That set is not empty because it contains the solution of the linear equation (\ref{2.12e}) with $v(0) = v_0$ and $V = 0$, obtained by a simplified version of Proposition 3.4. Clearly the fixed point $v$ satisfies the equation (\ref{1.21e}) and therefore belongs to $F(T, v_0)$.\par\nobreak \hfill $\sq$ \par

We finally translate the main result of Proposition 4.1 in terms of $u_c$ and we derive additional bounds and regularity properties for $u_c$.\\

\noi {\bf Proposition 4.2.} {\it Let $1 - \gamma /2 < \rho < n/2$, let $v_0 \in H^\rho$ and define $\varphi$ by (\ref{1.18e}). Then there exists $T > 0$ and there exists a unique solution $u_c \in {\cal C}([0, T], H^\rho )$ of the equation (\ref{1.11e}) such that $v$ defined by (\ref{1.12e}) satisfies the equation (\ref{1.21e}) with $v(0) = v_0$. Furthermore $u_c$ satisfies the estimate 
 \beq
\label{4.14e}
\parallel u_c(t); H^\rho  \parallel \ \leq  C\ a_0 \left ( 1 + a_0^2 \ t^{\gamma - 1} \right )^{1 + [\rho ]} 
\eeq

\noi for all $t \in (0, T]$, where $[ \rho ]$ is the integral part of $\rho$ and}
$$a_0 \ = \ \parallel v_0; H^\rho  \parallel \ .$$
\vskip 3 truemm

\noi {\bf Proof.} The first statement follows from Proposition 4.1, except for the continuity of $u_c$. We first prove the estimate (\ref{4.14e}). By Lemmas 2.1 and 2.2 we estimate
 \beq
\label{4.15e}
\parallel \omega^{\rho } \exp (- i \varphi ) v) \parallel \ \leq \ C \left ( 1 \ + \ \parallel \omega^{\rho } (\exp (-i \varphi ) - 1) \parallel_{n/\rho} \right ) \parallel \omega^{\rho } v\parallel \ . 
\eeq 

\noi It follows from Lemma 2.5 that for $0 < \rho < n/2$ 
\bea
\label{4.16e}
&&\parallel \omega^{\rho }  ( \exp (- i \varphi ) - 1 ) \parallel_{n/\rho} \ \leq \ C \parallel \exp (- i \varphi ) - 1; \dot{B}_{n/\rho , 2}^\rho   \parallel \nn \\
&&\leq \ C \parallel \varphi  ; \dot{B}_{n/\rho , 2}^\rho \parallel \left (  1 \ + \parallel \varphi ; \dot{B}_{\infty , \infty}^0  \parallel \right )^{[\rho ]} \nn \\
&&\leq \ C\parallel \omega^{n/2} \varphi \parallel \left ( 1 \ + \parallel \omega^{n/2} \varphi\parallel \right )^{[\rho ]}   \ .                                                                                                                                                                                                                                                                                                                                                                                                                                                                                                                                                                                                                                                                                                                                                                                                                                                                                                                                                                                                                                                                                                                                                                                                                                                                                                                                                                                                                                                                                                                                                                                                                                                                                                                                                                                                                                                                                                                                                                                                                                                                                                                         \eea
 
\noi Using (\ref{1.18e}) and Lemma 2.3, we estimate
\beq
\label{4.17e}
\parallel \omega^{n/2} \varphi \parallel  \  \leq \ C \ t^{\gamma - 1} \parallel \omega^{\gamma /2} v_0\parallel^2\ \leq \ C\ a_0^2 \ t^{\gamma - 1}  
\eeq

\noi which together with (\ref{4.15e}) and (\ref{4.16e}) implies (\ref{4.14e}).\par

It remains to prove the continuity in time of $u_c$ in $H^\rho$. For that purpose it suffices to show that the multiplication by $\exp (- i \varphi )$ is strongly continuous in $t$ as an operator in $H^\rho$. Now for fixed $v$ and $t_0$
\bea
\label{4.18e}
&&\parallel  \left ( \exp (- i \varphi (t)) - \exp (- i \varphi (t_0))\right )  v ; H^\rho  \parallel \nn \\
&&\leq \ C \left ( \parallel  \delta \varphi  \parallel_\infty  \ + \ \parallel \omega^{n/2} \delta \varphi \parallel  \left ( 1 \  + \ \parallel \omega^{n/2} \delta \varphi \parallel \right )^{[\rho ]} \right ) \parallel \exp (- i \varphi (t_0)) v ; H^\rho \parallel \nn \\   \eea

\noi where $\delta \varphi = \varphi (t) - \varphi (t_0)$ and it suffices to prove that $\varphi (t)$ is a continuous function of $t$ in $L^\infty \cap \dot{H}^{n/2}$. This follows immediately from estimates similar to (\ref{4.17e}).\par\nobreak \hfill $\sq$ \par
   
\mysection{Continuity with respect to initial data}
\hspace*{\parindent} 
In this section we prove that the map $v_0 \to v$ defined by Proposition 4.1 is continuous in the natural norms and that the map $v_0 \to u_c$ defined by Proposition 4.2 satisfies similar continuity properties. As mentioned in the introduction, the proof of those properties given here requires the additional condition $\rho > 3/4$, which is stronger than the condition $\rho > 1 - \gamma /2$ for $\gamma > 1/2$. We need to estimate the difference of two solutions of the linearized equation (\ref{1.23e}) corresponding to two functions $v_1$ and $v_2$ not necessarily satisfying the condition $v_1(0) = v_2(0)$. The following lemma is an extension of Lemma 4.1 where we drop that condition. Furthermore we assume for simplicity that $\rho ' = \rho$. \\

\noi {\bf  Lemma 5.1.} {\it Let $3 /4 < \rho < n/2$. Let $I = (0,T]$ and let $v_i$, $i = 1,2$, satisfy the assumptions of Proposition 3.4 with $v_i(0) = v_{0i} \in H^\rho$. Let $v'_i$, $i = 1,2$ be the solutions of the equation (\ref{1.23e}) with $v'_i (0) = v'_{0i} \in H^\rho$ obtained in Proposition 3.4. Let $0 \leq \sigma \leq \rho - 1/2$ with $\rho + \sigma > 1$. Define 
$$y \equiv y(t) = \ \parallel \omega^\sigma v_- (t)\parallel^2 \ + \  \parallel v_- (t)\parallel^2\ \approx \ \parallel  v_- (t) ; H^\sigma \parallel^2\ $$
$$y' \equiv y'(t) = \ \parallel \omega^\sigma v'_- (t)\parallel^2 \ + \  \parallel v'_- (t)\parallel^2\ \approx \ \parallel  v'_- (t) ; H^\sigma \parallel^2$$
\noi and $y_0 = y(0)$. Then the following estimate holds
\bea  
\label{5.1e}
&&\left | \partial_t y' \right | \leq C\ t^{\lambda_1 - 1} \left ( a^2 y' + aa{'}^2 y_0^{1/2} + aa' y{'}^{1/2} \Big ( y_0^{1/2} + y^{1/2} + \ t^{-1} \int_0^t dt' \ y^{1/2} (t') \Big ) \right ) \nn \\
&&+ \ C\ t^{2\lambda_0 - 2} \left ( a^4 y' + a^3 a' y{'}^{1/2} \Big ( y_0^{1/2} + t^{-\lambda_0} \int_0^t dt' \ t{'}^{\lambda_0 - 1} y^{1/2} (t')\Big ) \right )
 \eea	
\noi for all $t \in I$, where $\lambda_0$, $\lambda_1$ are defined by (\ref{2.10e}), and $a$, $a'$ are defined by (\ref{4.2e}).}\\

\noi {\bf Proof.} In the same way as in the proof of Lemma 4.1, we obtain from (\ref{1.23e})
\beq
\label{5.2e}
i \partial_t v'_- = L_2 \ v'_- + L_- \ v'_1
\eeq

\noi where $L_i = L(v_i)$,
\beq
\label{5.3e}
L_- = i \ s_- \cdot \nabla + (i/2) (\nabla \cdot s_-) + f_-
\eeq

\noi and $f$ is defined by (\ref{3.6e}). We estimate for $0 \leq \sigma \leq \rho$ 
\beq
\label{5.4e}
\partial_t \parallel \omega^{\sigma} v'_-\parallel^2 \ = 2 \ {\rm Im} \left ( <  \omega^\sigma v'_- ,  \omega^\sigma L_2 v'_-> \ + \  <  \omega^\sigma v'_- ,  \omega^\sigma L_-v'_1> \right ) \ .
\eeq
 
\noi We estimate the first scalar product in (\ref{5.4e}) as in the proof of Proposition 3.2 namely (see (\ref{3.13e}) (\ref{3.14e})).
\beq
\label{5.5e}
 \left | <  \omega^\sigma v'_- ,  \omega^\sigma L_2 v'_-> \right | \leq \ N(t)  \parallel \omega^{\sigma} v'_-\parallel^2 \ .
\eeq

\noi We next estimate the second scalar product in (\ref{5.4e}) and we estimate the contribution of the various terms of (\ref{5.3e}) successively. We first estimate the contribution of $s_-\cdot \nabla$, namely
\beq
\label{5.6e}
M =  \left | <  \omega^\sigma v'_- ,  \omega^\sigma s_-\cdot \nabla v'_1> \right | \ .
\eeq

\noi We consider separately the cases $\sigma \geq \rho - 1$ and $\sigma < \rho - 1$. For $\sigma \geq \rho - 1$, we estimate
\beq
\label{5.7e}
M \leq\ C   \parallel \omega^{2\sigma + 1 - \rho} \ v'_-\parallel \left ( \parallel s_- \parallel_{\infty} \ + \ \parallel \omega^{n/2} s_-\parallel \right )  \parallel \omega^\rho v'_1\parallel
\eeq

\noi by Lemma 2.2. We next try to estimate
\bea
\label{5.8e}
\parallel s_-\parallel_\infty \ + \ \parallel \omega^{n/2} s_-\parallel &\leq&  C \ t^{\gamma -1} \parallel \omega^{1 + \gamma - n/2 \pm \varepsilon} \chi_L |v_{0} |_-^2 \parallel \nn \\
&\leq&  C \ t^{\gamma - 1} \parallel \omega^{\rho_{2\pm}}  v_{0+} \parallel \  \parallel \omega^{\sigma_2}  v_{0-} \parallel 
\eea

\noi by Lemma 2.3, where $|v_0|_-^2 \equiv (|v_0|^2)_-$, with
\beq
\label{5.9e}
0 \leq \rho_{2\pm} \leq \rho \ , \ 0 \leq \sigma_2 \leq \sigma \ , \ \rho_{2\pm} + \sigma_2 = 1 + \gamma \pm \varepsilon \ ,
\eeq

\noi which implies $\sigma_2 \leq 1 + \gamma - \varepsilon$. We choose $\sigma_2 = \sigma \wedge (1 + \gamma - \varepsilon )$. \par

If $\sigma \geq 1 + \gamma - \varepsilon$, then $\sigma_2 = 1 + \gamma - \varepsilon$ and $\rho_{2\pm} = \varepsilon \pm \varepsilon \leq \rho$.

If $\sigma \leq 1 + \gamma - \varepsilon$, then $\sigma_2 = \sigma$ and $\rho_{2\pm} = 1 + \gamma -\sigma \pm \varepsilon$, which satisfies $\rho_{2\pm} \leq \rho$ for $\rho + \sigma > 1 + \gamma$. For $\rho + \sigma < 1 + \gamma$, we use the cut off $\chi_L$ to replace $\rho_{2\pm}$ by $\rho$, so that
\bea
\label{5.10e}
\parallel s_-\parallel_\infty \ + \ \parallel \omega^{n/2} s_-\parallel &\leq&  C \ t^{\gamma -1 - (1/2)[1 + \gamma - \rho - \sigma]_+} \parallel v_{0+} ; H^\rho \parallel \nn \\
&\times& \parallel \omega^{\sigma \wedge (1 + \gamma - \varepsilon )} v_{0-} \parallel  \ . 
\eea

\noi Now $\sigma \geq \rho - 1$ implies $1 + \gamma - \rho - \sigma \leq 2 + \gamma - 2\rho$. Substituting (\ref{5.10e}) into (\ref{5.7e}) then yields
\beq
\label{5.11e}
M \leq\ C \  t^{\lambda_1 - 1} aa'\parallel \omega^{2\sigma + 1 - \rho}  v'_-\parallel \ \parallel  \omega^{\sigma \wedge (1 + \gamma - \varepsilon )}  v_{0-}\parallel \  .
\eeq

\noi The first norm in (\ref{5.11e}) is controlled by the $H^\rho$ norm of $v'_-$ only if $\sigma \leq \rho - 1/2$, and that is the origin of that condition, which is otherwise not used to derive (\ref{5.11e}).\par

We next estimate $M$ in the case where $\sigma < \rho - 1$ which we rewrite as $\sigma \leq \rho - 1 - \varepsilon$. We estimate
\beq
\label{5.12e}
M \leq\ C  \parallel \omega^{\sigma} v'_-\parallel \ \parallel  \omega^{n/2 - \rho_1 + 1 + \sigma} s_{-}\parallel \  \parallel  \omega^{\rho_1}  v'_{1}\parallel
\eeq 

\noi by Lemma 2.3 for
\beq
\label{5.13e}
0 \leq \sigma \leq \rho_1 - 1 - \varepsilon \ .
\eeq

\noi We next try to estimate 
\bea
\label{5.14e}
&&\parallel \omega^{n/2 - \rho_1 + 1 + \sigma} s_-\parallel \ = \ C \  t^{\gamma - 1} \parallel  \omega^{2 + \gamma - \rho_1 + \sigma - n/2} \ \chi_L |v_0|_-^2\parallel \nn \\
&&\leq C\ t^{\gamma -1}  \parallel \omega^{\rho_2}  v_{0+}\parallel \ \parallel \omega^{\sigma_2}  v_{0-}\parallel
\eea

\noi by Lemma 2.3 again, with
\beq
\label{5.15e}
0 \leq \rho_{1}, \rho_2 \leq \rho \ , \ 0 \leq \sigma_2 \leq \sigma \ , \ \rho_{2} + \sigma_2 = 2 + \gamma - \rho_1 + \sigma 
\eeq

\noi which together with (\ref{5.13e}), implies $\sigma_2 \leq 1 + \gamma - \varepsilon$. We choose again $\sigma_2 = \sigma \wedge (1 + \gamma - \varepsilon )$. \par

If $\sigma \geq 1 + \gamma - \varepsilon$, then $\sigma_2 = 1 + \gamma - \varepsilon$ and $\rho_{1} + \rho_2 = 1 + \sigma + \varepsilon$. We choose $\rho_2 = 0$ and $\rho_1 = 1 + \sigma + \varepsilon$, so that $\rho_2 < \rho_1 \leq \rho$.\par

If $\sigma \leq 1 + \gamma - \varepsilon$, then $\sigma_2 = \sigma$ and $\rho_{1} + \rho_2 = 2 + \gamma$. We choose
\beq
\label{5.16e}
\left \{ \begin{array}{l} \rho_1 = (1 + \gamma /2) \vee (1 + \sigma + \varepsilon ) \\ \\ \rho_2 = (1 + \gamma /2) \wedge (1 + \gamma -\sigma - \varepsilon ) \end{array} \right .
\eeq

\noi which implies $\rho_2 \leq \rho_1 \leq \rho$ and therefore ensures (\ref{5.14e}), except in the case where $\sigma < \gamma /2 - \varepsilon$ and $\rho < 1 + \gamma /2$. In that case, we use the cut off $\chi_L$ to replace $\rho_1$ and $\rho_2$ by $\rho$, so that 
\beq
\label{5.17e}
\parallel \omega^{n/2 - \rho + 1 + \sigma} s_-\parallel \ \leq \ C \  t^{\lambda_1 - 1} \parallel  \omega^{\rho}  v_{0+}\parallel  \ \parallel \omega^{\sigma} v_{0-}\parallel \ .
\eeq

\noi Substituting (\ref{5.14e}) or (\ref{5.17e}) into (\ref{5.12e}), together with (\ref{5.11e}), yields the estimate
\beq
\label{5.18e}
M \leq\ C \  t^{\lambda_1 - 1} aa'\parallel \omega^{\sigma \vee (2 \sigma + 1 - \rho )}  v'_-\parallel \ \parallel  \omega^{\sigma \wedge (1 + \gamma - \varepsilon )} v_{0-}\parallel 
\eeq

\noi for $0 \leq \sigma < \rho$.\par

We next turn to the remaining terms from (\ref{5.3e}). They take the form
$$<\omega^\sigma v'_-, \omega^\sigma F\ v'_1>$$

\noi for some function $F$ and are estimated as
\beq  
\label{5.19e}
\left | <\omega^{\sigma} v'_- , \omega^{\sigma} F v'_1 >\right | \ \leq \ C \parallel \omega^{\sigma } v'_- \parallel \ \parallel \omega^{n/2 - \rho_1 + \sigma} F \parallel \  \parallel \omega^{\rho_1} v'_1 \parallel
\eeq 

\noi by Lemma 2.3, for some $\rho_1$, $0 \leq \rho_1 \leq \rho$, with
\beq
\label{5.20e}
0 \leq \sigma \leq \rho_1 - \varepsilon \ .
\eeq

\noi We now estimate the middle norm in the RHS of (\ref{5.19e}) for the relevant choices of $F$ and suitable choices of $\rho_1$. \par

We first consider the contribution of $\nabla\cdot s_-$. We try to estimate
\bea
\label{5.21e}
\parallel \omega^{n/2 - \rho_1 + \sigma} \ \nabla \cdot s_-\parallel &\leq& C \ t^{\gamma - 1}\parallel  \omega^{2 + \gamma - \rho_1 + \sigma - n/2}  \ \chi_L |v_0|_-^2\parallel \nn \\  
&\leq& C\ t^{\gamma - 1}\parallel  \omega^{\rho_2}   v_{0+}\parallel \ \parallel  \omega^{\sigma_2}   v_{0-}\parallel
\eea

\noi with $\rho_1$, $\rho_2$, $\sigma_2$ satisfying (\ref{5.20e}) (\ref{5.15e}), which imply $\sigma_2 \leq 2 + \gamma - \varepsilon$. We choose $\sigma_2 = \sigma \wedge (2 + \gamma - \varepsilon )$. \par

If $\sigma \geq 2 + \gamma - \varepsilon$, we choose $\rho_1 = \sigma + \varepsilon$ and $\rho_2 = 0$. \par

If $\sigma \leq 2 + \gamma - \varepsilon$, so that $\sigma_2 = \sigma$, we choose
\beq  
\label{5.22e}
\left \{ \begin{array}{l} \rho_1 = (1 + \gamma /2 ) \vee ( \sigma  + \varepsilon ) \\ \\ \rho_2 = (1 + \gamma /2 ) \wedge ( 2 + \gamma - \sigma  - \varepsilon ) \end{array} \right .
\eeq 

\noi which ensures (\ref{5.20e}) (\ref{5.15e}) and (\ref{5.21e}) for $\rho \geq 1 + \gamma /2$. For $\rho < 1 + \gamma /2$, we use the cut off $\chi_L$ to replace $\rho_1$ and $\rho_2$ by $\rho$, so that finally
\beq
\label{5.23e}
\parallel \omega^{n/2 - \rho_1 + \sigma} \ \nabla \cdot s_-\parallel \ \leq \ C \  t^{\lambda_1 - 1} a\parallel  \omega^{\sigma \wedge (2 + \gamma - \varepsilon)}   v_{0-}\parallel  \ .
\eeq

We next consider the contribution of the various terms of $f_-$, where $f$ is defined by (\ref{3.6e}). We first consider the contribution of $|s|_-^2 (\equiv (|s|^2)_-)$. We estimate
\beq
\label{5.24e}
\parallel \omega^{n/2 - \rho_1 + \sigma} |s|_-^2\parallel \ \leq \ C \parallel  \omega^{n/2 \pm 0}  s_{+}\parallel \  \parallel  \omega^{n/2 - \rho_1 +\sigma}  s_-\parallel
\eeq

\noi by Lemma 2.2. The first norm in the RHS is estimated by (\ref{2.9e}). We try to estimate the second norm by
\bea  
\label{5.25e}
&&\parallel \omega^{n/2 - \rho_1 + \sigma} s_- \parallel \ \leq \  C\ t^{\gamma - 1}  \parallel \omega^{1 + \gamma - \rho_1 + \sigma - n/2} \ \chi_L |v_0|_-^2\parallel \nn \\
&&\leq C\ t^{\gamma - 1} \parallel \omega^{\rho_2}  v_{0+} \parallel \ \parallel \omega^{\sigma_2}  v_{0-} \parallel
 \eea	

\noi with $\rho_1$, $\rho_2$, $\sigma_2$ satisfying (\ref{5.20e}) and
\beq
\label{5.26e}
0 \leq \rho_{1}, \rho_2 \leq \rho \ , \ 0 \leq \sigma_2 \leq \sigma \ , \ \rho_{2} + \sigma_2 = 1 + \gamma - \rho_1 + \sigma 
\eeq

\noi which imply $\sigma_2 \leq 1 + \gamma - \varepsilon$. We choose $\sigma_2 = \sigma \wedge (1 + \gamma - \varepsilon )$. \par

If $\sigma \geq 1 + \gamma - \varepsilon$, we choose $\rho_1 = \sigma + \varepsilon$ and $\rho_2 = 0$. \par

If $\sigma \leq 1 + \gamma - \varepsilon$, so that $\sigma_2 = \sigma$, we choose
\beq  
\label{5.27e}
\left \{ \begin{array}{l} \rho_1 = (1 + \gamma ) /2 \vee ( \sigma  + \varepsilon ) \\ \\ \rho_2 = (1 + \gamma )/2  \wedge ( 1 + \gamma - \sigma  - \varepsilon ) \end{array} \right .
\eeq

\noi which ensures (\ref{5.20e}) (\ref{5.26e}) and (\ref{5.25e}) for $\rho \geq (1 + \gamma )/2$. For $\rho < (1 + \gamma )/2$, we use the cut-off $\chi_L$ to replace $\rho_1$ and $\rho_2$ by $\rho$ so that finally
\beq
\label{5.28e}
\parallel \omega^{n/2 - \rho_1 + \sigma} s_-\parallel \ \leq \ C \ t^{\lambda_0 - 1} \ a\parallel  \omega^{\sigma \wedge (1 + \gamma - \varepsilon )}   v_{0-}\parallel 
\eeq

\noi and 
\beq
\label{5.29e}
\parallel \omega^{n/2 - \rho_1 + \sigma} |s|_-^2\parallel \ \leq \ C \ t^{2\lambda_0 - 2} \ a^3\parallel  \omega^{\sigma \wedge (1 + \gamma - \varepsilon )} v_{0-}\parallel \ .
\eeq

\noi We next estimate the contribution of $g_{S-}$ to $f_-$ by
\bea  
\label{5.30e}
&&t^{\gamma -2}\parallel \omega^{n/2 - \rho + \sigma} g_{S-} \parallel \ = \ C\ t^{\gamma - 2}  \parallel \omega^{\gamma - \rho + \sigma - n/2} \ \chi_S |v|_-^2\parallel \nn \\
&&\leq \ C\ t^{\gamma - 2 + \rho - \gamma /2} \parallel \omega^{\rho} v_{+} \parallel \ \parallel \omega^{\sigma} v_{-} \parallel\nn \\
&&\leq \ C\ t^{\lambda_1 - 1} a \parallel \omega^{\sigma} v_{-} \parallel
 \eea	

\noi for $\rho_1 = \rho \geq \gamma /2$. \par

Similarly
\beq
\label{5.31e}
\parallel \omega^{n/2 - \rho + \sigma} \left ( t^{\gamma - 2} g_L (v_0)_- - \partial_t \varphi_-\right ) \parallel \ \leq \ C \ t^{\lambda_1 - 1} a\parallel  \omega^{\sigma}  v_{0-}\parallel \ .
\eeq

We next consider the contribution of $(g_L (v) - g_L (v_0))_-$ to $f_-$. We want to estimate
\bea  
\label{5.32e}
J &=&t^{\gamma -2}\parallel \omega^{n/2 - \rho_1 + \sigma} \left ( g_{L}(v) - g_L (v_0) \right )_- \parallel  \nn \\
&=&t^{\gamma - 2} \parallel \omega^{\gamma - \rho_1 + \sigma - n/2} \ \chi_{L} \left (  |v|_{-}^2 - |v_0|_-^2 \right ) \parallel \ .
 \eea	

\noi From the conservation law (\ref{2.16e}), we obtain
\beq
\label{5.33e}
 |v|_-^2 -  |v_0|_-^2 =  V_{1-} +  V_{2-} +  V_{3-}
\eeq

\noi where
\bea
\label{5.34e}
&&V_{1-} = - \int_0^t dt'\ {\rm Im} \left (  \overline{v}_+ \Delta v_- +  \overline{v}_- \Delta v_+ \right ) (t') \\
&&V_{2-} = \nabla \cdot \int_0^t dt'  \left (  s_+ |v|_-^2  \right ) (t') \\
\label{5.35e}
&&V_{3-} = \nabla \cdot \int_0^t dt'  \left (  s_- |v|_+^2  \right ) (t') \ .
\label{5.36e}
\eea

\noi We first consider the contribution of $V_{1-}$. By the same estimates as in Lemma 2.6, we obtain
\beq
\label{5.37e}
\parallel \omega^{\rho_2 + \sigma_2 - 2 - n/2} \ V_{1-} \parallel \ \leq \ C \ \int_0^t dt' \parallel  \omega^{\rho_2}   v_{+}\parallel \ \parallel  \omega^{\sigma_2}   v_{-}\parallel (t')
\eeq

\noi for $0 \leq \rho_2$, $\sigma_2 < n/2$, $1 < \rho_2 + \sigma_2 \leq 2 + n/2$. We try to estimate
\beq  
\label{5.38e}
t^{\gamma -2}\parallel \omega^{\gamma - \rho_1 + \sigma - n/2} \ \chi_L \ V_{1-}  \parallel \ \leq \ C\ t^{\gamma - 2}  \int_0^t dt' \left ( \parallel \omega^{\rho_2} v_+\parallel \ \parallel  \omega^{\sigma_2} v_-\parallel \right ) (t')
 \eeq	

\noi with $\rho_1$, $\rho_2$, $\sigma_2$ satisfying (\ref{5.20e})  (\ref{5.15e}) and $\rho_2 + \sigma_2 > 1$. We proceed as for the estimate of the contribution of $\nabla\cdot s_-$ (see (\ref{5.21e})  (\ref{5.23e})) and we obtain finally
\beq  
\label{5.39e}
t^{\gamma -2}\parallel \omega^{\gamma - \rho_1 + \sigma - n/2} \ \chi_L \ V_{1-}  \parallel \ \leq \ C\ t^{\lambda_1 - 2}  a \int_0^t dt'  \parallel \omega^{\sigma \wedge (2 + \gamma - \varepsilon)} v_-(t') \parallel \ . 
 \eeq	

\noi The condition $\rho_2 + \sigma_2 > 1$ is ensured provided $\rho + \sigma > 1$.\par

We shall also need an estimate of $\parallel v'_-\parallel$. In all terms but $V_{1-}$ this is done by taking $\sigma = 0$ in the available estimate. Doing so for $V_{1-}$ would require the stronger condition $\rho > 1$. We shall instead estimate the corresponding norm of $V_{1-}$ in terms of the norm $\parallel \omega^{\sigma_2} v_-\parallel$ for the previous choice of $\sigma_2$, thereby obtaining
\beq  
\label{5.40e}
t^{\gamma -2}\parallel \omega^{\gamma - \rho_1 - n/2} \ \chi_L \ V_{1-}  \parallel \ \leq \ C\ t^{\lambda_1 - 2}  a \int_0^t dt'  \parallel \omega^{\sigma \wedge (2 + \gamma - \varepsilon)}  v_-(t')\parallel 
 \eeq	

\noi by a similar computation. \par

We next consider the contribution of $V_{2-}$. We try to estimate  
\bea  
\label{5.41e}
&&t^{\gamma -2}\parallel \omega^{\gamma - \rho_1 + \sigma - n/2} \ \chi_L \ V_{2-}  \parallel \nn \\
&&\leq \  t^{\gamma - 2}  \int_0^t dt'  \parallel \omega^{1 + \gamma - \rho_1 + \sigma - n/2} \chi_L(t) \left ( s_+ |v|_-^2 \right ) (t')\parallel \nn \\ 
&&\leq \  C\ t^{\gamma - 2}  \int_0^t dt'  \left ( \left ( \parallel s_+ \parallel_\infty \ + \ \parallel \omega^{n/2} s_+ \parallel \right )  \parallel \omega^{\rho_2}  v_+\parallel \ \parallel \omega^{\sigma_2}  v_-\parallel \right ) (t')\nn \\ 
 \eea	

\noi with $\rho_1$, $\rho_2$, $\sigma_2$ satisfying (\ref{5.20e}) (\ref{5.26e}). We estimate the norms of $s_+$ by (\ref{2.9e}) and we proceed for the remaining norms as for the estimate of $s_-$ (see (\ref{5.25e})-(\ref{5.28e})), thereby obtaining finally
\bea  
\label{5.42e}
&&t^{\gamma -2}\parallel \omega^{\gamma - \rho_1 + \sigma - n/2} \chi_L V_{2-}  \parallel \nn \\
&& \leq \ C\ t^{\lambda_0 - 2}  \ a^3 \int_0^t dt' \ t{'}^{\lambda_0 - 1} \parallel \omega^{\sigma \wedge (1 + \gamma - \varepsilon)}  v_-(t')\parallel \ . 
 \eea	

\noi We next consider the contribution of $V_{3-}$. We estimate
\bea  
\label{5.43e}
&&t^{\gamma -2}\parallel \omega^{\gamma - \rho_1 + \sigma - n/2} \ \chi_L \ V_{3-}  \parallel \nn \\
&&\leq \ t^{\gamma - 2}  \int_0^t dt' \parallel \omega^{1 + \gamma - \rho_1 + \sigma - n/2} \ \chi_L(t) \left ( s_-|v|_+^2 \right ) (t')\parallel \ . 
\eea	

\noi For $\rho \geq (1 + \gamma )/2$, we estimate the last norm by
\bea  
\label{5.44e}
\parallel \ \cdot \ \parallel &\leq& C\ t{'}^{\lambda_0 -1}\parallel \omega^{1 + \gamma - n/2} |v|_+^2 \parallel \  \parallel \omega^{n/2 - \rho_1 + \sigma}  s_-\parallel \nn \\
&\leq& C \ t{'}^{\lambda_0 - 1} \parallel \omega^{(1 + \gamma )/2}  v_{+} \parallel^2 \ \parallel \omega^{\rho_2}  v_{0+} \parallel
\  \parallel \omega^{\sigma \wedge (1 + \gamma - \varepsilon )}  v_{0-} \parallel
 \eea	

\noi by Lemma 2.3 and by (\ref{5.25e}) with the choice (\ref{5.27e}) of $\rho_1$, $\rho_2$. For $\rho < (1 + \gamma )/2$, we take $\rho_1 = \rho$ and we use the cut off $\chi_L$ twice to estimate
\bea  
\label{5.45e}
&&t^{\gamma -2}\parallel \omega^{1 + \gamma - \rho + \sigma - n/2} \ \chi_L (t) \left ( s_- |v|_+^2 \right ) (t')  \parallel \nn \\
&& \leq \ C \ t^{\lambda_0 - 2}  \parallel \omega^{\rho + \sigma - n/2} \left ( s_-|v|_+^2 \right ) (t')\parallel \nn \\
&& \leq \ C\ t^{\lambda_0 - 2}  \parallel \omega^{2\rho - n/2} |v|_+^2 \parallel \ \parallel \omega^{\sigma - \rho + n/2} s_-(t') \parallel \nn \\
&&\leq \ C\ t^{\lambda_0 - 2}  \ t{'}^{\lambda_0 - 1} \parallel \omega^{\rho} v_+\parallel^2  \ \parallel \omega^{\rho + \sigma - n/2} \ |v_0|_-^2 \parallel \nn \\
&&\leq \ C\ t^{\lambda_0 - 2}  \ t{'}^{\lambda_0 - 1} \parallel \omega^{\rho} v_+\parallel^2  \ \parallel \omega^{\rho} v_{0+} \parallel  \ \parallel \omega^{\sigma} v_{0-} \parallel
\eea	

\noi by repeated use of Lemma 2.3. Substituting (\ref{5.44e}) or (\ref{5.45e}) into (\ref{5.43e}) and integrating over time yields
\beq  
\label{5.46e}
t^{\gamma -2}\parallel \omega^{\gamma - \rho_1 + \sigma - n/2} \ \chi_L \ V_{3-}  \parallel \ \leq \ C \ t^{2\lambda_0 - 2}  \ a^3 \parallel \omega^{\sigma \wedge (1 + \gamma - \varepsilon )}  v_{0-}\parallel \ . 
\eeq	

\noi We substitute (\ref{5.23e}) (\ref{5.29e}) (\ref{5.30e}) (\ref{5.31e}) (\ref{5.40e}) (\ref{5.42e}) (\ref{5.46e}) into (\ref{5.19e}) and substitute the result as well as (\ref{5.5e}) and (\ref{5.18e}) into (\ref{5.4e}), thereby obtaining
\bea  
\label{5.47e}
&&\left | \partial_t  \parallel \omega^{\sigma} v'_{-}(t) \parallel^2\right | \leq \ N(t)  \parallel \omega^{\sigma} v'_{-} \parallel^2 \nn \\
&&+ \ C\ aa{'}\ t^{\lambda_1 - 1} \Big ( \parallel \omega^{\sigma \vee (2\sigma + 1 - \rho )} \ v'_{-}\parallel \ \parallel \omega^{\sigma \wedge (1 + \gamma )_-} \ v_{0-}\parallel\nn \\
&&+\  \parallel \omega^{\sigma} v'_{-} \parallel \Big ( \parallel \omega^{\sigma \wedge (2 + \gamma )_-} \ v_{0-}\parallel \ + \ \parallel \omega^{\sigma} v_{-} \parallel \ + \ \parallel \omega^{\sigma} v_{-} \parallel \nn \\
 &&+ \ t^{-1} \int_0^t dt' \parallel \omega^{\sigma \wedge (2 + \gamma )_-} \ v_{-}(t')\parallel \Big ) \Big ) \nn \\
&&+ \ C\ a^3a{'}\ t^{2\lambda_0 - 2} \parallel \omega^{\sigma} v'_{-}\parallel \Big (  \parallel \omega^{\sigma \wedge (1 + \gamma )_-} \ v_{0-}\parallel\nn \\
&&+ \  t^{-\lambda_0} \int_0^t dt' \ t{'}^{\lambda_0 - 1} \parallel \omega^{\sigma \wedge (1 + \gamma )_-} \ v_{-}(t')\parallel \Big ) 
 \eea	

\noi where $(j + \gamma )_- = (j + \gamma - \varepsilon )$ for $j = 1,2$. Together with the similar estimate for $\sigma = 0$ (see however the discussion after (\ref{5.39e})), with (\ref{3.14e}) and the fact that the second norm of $v'_-$ in the RHS is bounded by $a'$ for $\sigma \leq \rho - 1/2$, this yields (\ref{5.1e}). \par\nobreak \hfill $\sq$\\

\noi {\bf Remark 5.1.} All the estimates leading to (\ref{5.47e}) hold for $0 \leq \sigma < \rho$ except for the estimate of $V_{1-}$ which requires $\rho + \sigma > 1$. On the other hand the estimate (\ref{5.18e}) coming from the transport term $s_- \cdot \nabla$ can be used only for $\sigma \leq \rho - 1/2$. Those two conditions force the restriction $\rho > 3/4$. \\

We can now state the continuity properties of the map $v_0 \to v$ defined in Proposition 4.1.\\

\noi {\bf Proposition 5.1.} {\it Let $3/4 < \rho < n/2$. Let $R > 0$ and let $T$ be defined by (\ref{4.12e}). Let $B_0 (R/2)$ be the ball of radius $R/2$ in $H^\rho$.

(1) Let $1 - \rho < \sigma < \rho$ and $\sigma \geq 0$. Then the map $v_0 \to v$ defined by Proposition 4.1 is continuous from $H^{\sigma}$ to $L^\infty ((0, T], H^\sigma )$ uniformly for $v_0 \in B_0 (R/2)$. Furthermore, for $1 - \rho < \sigma \leq \rho - 1/2$ and for two solutions $v_i$, $i = 1,2$ of the equation (\ref{1.21e}) with $v_i(0) \equiv v_{i0} \in B_0 (R/2)$ as obtained in that proposition, the following estimate holds for all $t\in (0, T]$
 \beq
\label{5.48e}
\parallel v_-(t) ;  H^{\sigma} \parallel^2  \ \leq \left ( 1 + C\ t^{\lambda}  \exp \left ( C t^{\lambda}\right ) \right )  \left ( y_0 + C\ t^{\lambda} \left (y_0 + y_0^{1/2}\right  ) \right )
\eeq 

\noi where $\lambda = \lambda_1 \wedge (2 \lambda_0 - 1)$ and $y_0 = \parallel v_{0-} ; H^{\sigma} \parallel^2$.\par

For $\rho - 1/2 < \sigma < \rho$, a similar estimate can be obtained by interpolation between (\ref{5.48e}) with $\sigma = \rho - 1/2$ and boundedness in $H^\rho$.\par

(2) The map $v_0 \to v$ defined by Proposition 4.1 is (pointwise) continuous from $H^\rho$ to $L^\infty ((0, T], H^\rho )$ for $v_0 \in B_0 (R/2)$}.\\

\noi {\bf Proof.} Part (1). Let $v_i$, $i = 1,2$, be two solutions of the equation (\ref{1.21e}) as defined above. Then (\ref{5.1e}) with $y' = y$ and $a' = a$ yields  
\bea  
\label{5.49e}
\left | \partial_t y \right | &\leq&  C\ t^{\lambda - 1} \Big ( y + y^{1/2} \ y^{1/2}_0 + y^{1/2}_0 + y^{1/2} \ t^{-1} \int_0^t dt' y^{1/2}(t') \nn \\
&&+ \ y^{1/2} \ t^{-\lambda_0} \int_0^t dt'  t{'}^{\lambda_0 - 1} y^{1/2}(t') \Big ) \  .
 \eea	

\noi (The constant $C$ depends on $a$ through a factor $(a^2 + a^4)$). Using the inequalities
\begin{eqnarray*}  
&&t^{-1} \int_0^t dt' y^{1/2}(t') \leq \left (  t^{-1} \int_0^t dt' y(t') \right )^{1/2} \\
&&t^{-\lambda_0} \int_0^t dt'  t{'}^{\lambda_0 - 1} y^{1/2}(t') \leq \left ( \left ( 2 \lambda_0 - 1\right )^{-1}  t^{-1} \int_0^t dt' y(t') \right )^{1/2} 
 \end{eqnarray*}	

\noi yields 
\beq  
\label{5.50e}
\left | \partial_t y \right | \leq  C\ t^{\lambda - 1} \left ( \eta + y + \lambda^{-1} \ t^{-1} \int_0^t dt'  y(t') \right ) 
 \eeq

\noi where $\eta = y_0 + y_0^{1/2}$. Integration of (\ref{5.50e}) over time yields 
\beq
\label{5.51e}
y \leq y_0 + C\left ( t^\lambda \ \eta + z \right )
\eeq

\noi where the new constant $C$ now depends on $\lambda$, we have used the fact that $0 < \lambda < 1$, and
\beq
\label{5.52e}
z = \int_0^t dt' \  t{'}^{\lambda - 1} y(t')
\eeq

\noi so that
\beq  
\label{5.53e}
\partial_t z  \leq  t^{\lambda - 1} \ y_0 + C \left ( t^{2\lambda -1} \ \eta + t^{\lambda -1} z \right ) \ .
 \eeq

\noi Integrating (\ref{5.53e}) and substituting the result into (\ref{5.51e}) yields (\ref{5.48e}) from which Part (1) follows. \\

\noi \underline{Part (2)}. The proof is identical with that of Part (2) of Proposition 5.1 in \cite{4r}.\par \nobreak \hfill $\sq$ \\

We finally prove the continuity of the map $v_0 \to u_c$ that follows from Proposition 5.1. \\

\noi {\bf Proposition 5.2.} {\it Let $3/4 < \rho < n/2$. Let $R > 0$ and let $T$ be defined by (\ref{4.12e}). Let $B_0 (R/2)$ be the ball of radius $R/2$ in $H^\rho$. Then the map $v_0 \to u_c$ defined by Proposition 4.2 is continuous from $H^\rho$ to $L^\infty ([t_1, T], H^\rho )$ for $v_0 \in B_0 (R/2)$ and for any $t_1$, $0 < t_1 < T$.}\\

\noi {\bf Proof.} By Proposition 5.1, part (2) and (\ref{1.12e}), it suffices to prove that the multiplication by $\exp (-i \varphi )$ is strongly continuous from $v_0 \in H^\rho$ as an operator in $L^\infty ([t_1, T], H^\rho )$. Let $v_{oi} \in H^\rho$, $i = 1,2$ and let $\varphi_i$ be the associated phases defined by (\ref{1.18e}). For fixed $v \in L^\infty ([t_1, T], H^\rho )$ and for all $t\in [t_1, T]$, we estimate 
\begin{eqnarray*}
&&\parallel \left ( \exp (-i \varphi_2 ) - \exp (-i \varphi_1 )\right ) v ; H^\rho \parallel \ \leq \\
&& C\left ( \parallel \delta \varphi \parallel_\infty \ + \ \parallel \omega^{n/2} \delta \varphi \parallel \left ( 1 + \parallel \omega^{n/2} \delta \varphi\parallel \right )^{[\rho ]}\right ) \parallel  \exp (-i \varphi_1 ) v ; H^\rho\parallel
\end{eqnarray*}

\noi where $\delta \varphi = \varphi_2 - \varphi_1$ and it suffices to prove that $\varphi$ is a continuous function of $v_0$ in $L^\infty ( [t_1, T], L^\infty \cap \dot{H}^{n/2})$. This follows from the fact that $\varphi$ is quadratic in $v_0$ and from estimates similar to (\ref{4.17e}). \par\nobreak \hfill $\sq$

\mysection*{Appendix A1}
\hspace*{\parindent}
In this appendix, we prove a slightly more general version of Lemma 2.4 where we drop the requirement that the estimating spaces are homogeneous under dilations. This extension would be useful to treat the main problem of this paper in spaces $H^\rho$ with $\rho \geq n/2$. \par

We introduce the usual notation for the standard Paley-Littlewood decomposition. We use the notation $\widehat{f}$ for the Fourier transform of $f$. Let $\widehat{\psi}_0 \in {\cal C}_0^\infty ({I\hskip - 1 truemm R}^n, {I\hskip - 1 truemm R}^+)$, $\widehat{\psi}_0 (\xi ) = 1$ for $|\xi | \leq 1$, $\widehat{\psi}_0 (\xi ) = 0$ for $|\xi | \geq 2$. We define $\widehat{\varphi}_0 (\xi ) = \widehat{\psi}_0 (\xi ) - \widehat{\psi}_0 (2\xi )$, $\widehat{\psi}_j (\xi ) = \widehat{\psi}_0 (2^{-j} \xi )$ and $\widehat{\varphi}_j (\xi ) = \widehat{\varphi}_0 (2^{-j}\xi )$ for all $j \in \Z$. For any positive integer $\nu$, we define
$$\widetilde{\varphi}_j^{(\nu )} = \sum_{|j - k| \leq \nu} \varphi_k\ .$$

\noi The superscript $\nu$ will be omitted for $\nu = 1$. For any $u \in {\cal S}'$, we define $u_j = \varphi_j \star u$, $\widetilde{u}_j^{(\nu )} = \widetilde{\varphi}_j^{(\nu )} \star u$ and $S_j(u) = \psi_j \star u$. We shall repeatedly use the estimates
$$\parallel \omega^{\lambda}  \widetilde{u}_j^{(\nu )}\parallel_r \ \leq \ \parallel \omega^{\lambda}  \widetilde{\varphi}_j^{(\nu + 1)} \parallel_1 \ \parallel \widetilde{u}_j^{(\nu )}\parallel_r \ = \ 2^{\lambda j}  \parallel \omega^{\lambda}  \widetilde{\varphi}_0^{(\nu + 1)} \parallel_1 \ \parallel \widetilde{u}_j^{(\nu )}\parallel_r \eqno({\rm A1.1})$$
\noi which holds for all $\lambda \in {I\hskip - 1 truemm R}$ and $1 \leq r \leq \infty$, and
$$\parallel \omega^{\lambda}  S_j (u)\parallel_r \ \leq \ \parallel \omega^{\lambda}  \psi_j \parallel_1 \ \parallel u\parallel_r \ = \ 2^{\lambda j}  \parallel \omega^{\lambda}  \psi_0\parallel_1 \ \parallel u\parallel_r \eqno({\rm A1.2})$$
 
\noi which holds for all $\lambda \geq 0$ and $1 \leq r \leq \infty$.\par

We shall use the following elementary lemma.\\

\noi {\bf Lemma A1.1.} {\it Let
$$M = \int d\xi \ d \eta \ \widehat{f}(\xi , \eta ) \ \widehat{u}(\xi ) \ \widehat{v} (\eta ) \ \widehat{m} (\xi - \eta )\ .$$

\noi Then
$$|M| \ \leq \ \parallel f\parallel_1 \ \parallel u\parallel_{r_1}  \ \parallel v\parallel_{r_2}\ \parallel m\parallel_{r_0} \eqno({\rm A1.3})$$

\noi for $1 \leq r_i \leq \infty$, $\displaystyle{\sum} 1/r_i = 1$.}\\

\noi {\bf Proof.} By the definition of the Fourier transform
$$M = \int dx\ dy \ f(x,y) \int dz\ u (-x-z) \ v(-y+z) \ m(z)$$

\noi so that
$$|M| \leq \int dx\ dy |f(x,y)| \ \mathrel{\mathop {\rm Sup}_{x,y}}\ \left | \int dz\ u(-x-z) v(-y+z) m(z)\right |$$

\noi from which (A1.3) follows by the H\"older inequality. \par\nobreak \hfill $\sq$\par

We want to estimate $<P_1u, [\omega^\lambda, m]P_2 v>$ which up to inessential complex conjugation reduces to
$$M = \int d\xi\ d\eta \ (P_1 \widehat{u}) (\xi )\  (P_2 \widehat{v})(\eta ) \ \widehat{m} (\xi - \eta ) \left ( |\xi |^\lambda - |\eta |^\lambda \right )$$
$$= \sum_{j,k,\ell} \int d\xi\ d\eta (P_1 \widehat{u}_j) (\xi ) \ (P_2 \widehat{v}_k)(\eta ) \ \widehat{m}_\ell (\xi - \eta ) \left ( |\xi |^\lambda - |\eta |^\lambda \right ) \eqno({\rm A1.4})$$

\noi by introducing the Paley-Littlewood decompositions of $u$, $v$, $m$. \par

The basic estimate is the following lemma.\\

\noi {\bf Lemma A1.2.} {\it Let $\lambda > 0$ and let $P_1$, $P_2$ satisfy the assumptions of Lemma 2.4. Then $M$ can be decomposed as a sum
$$M = M_1 + M_2 + M_3 + M_4\eqno({\rm A1.5})$$

\noi where the $M_i$'s satisfy the following estimates
$$|M_1| + |M_3| \ \leq \ C\parallel u; \dot{B}_{r_1, 2}^{\theta \mu} \parallel \ \parallel v \parallel_{r_2} \ \parallel m ; \dot{B}_{r_0, 2}^{(1 - \theta )\mu } \parallel\eqno({\rm A1.6})$$
$$|M_2| + |M_3| \ \leq \ C\parallel u \parallel_{r_1} \ \parallel v ; \dot{B}_{r_2, 2}^{\theta \mu}\parallel \ \parallel m ; \dot{B}_{r_0, 2}^{(1 - \theta )\mu} \parallel\eqno({\rm A1.7})$$
$$|M_4| \ \leq \ C\parallel u; \dot{B}_{r_1, 2}^{\theta (\mu - \nu)} \parallel \ \parallel v ; \dot{B}_{r_2, 2}^{(1 - \theta )(\mu - \nu)}\parallel \ \parallel \omega^{\nu - 1}\nabla m \parallel_{r_0} \eqno({\rm A1.8})$$

\noi where $\mu = \lambda + \alpha_1 + \alpha_2$, $0 \leq \nu \leq 1$, $1 \leq r_i \leq \infty$, $\displaystyle{\sum} 1/r_i = 1$ and $\theta \in {I\hskip - 1 truemm R}$. The parameters $r_i$ and $\theta$ can be chosen independently in the estimates (A1.6)-(A1.8).}\\

\noi {\bf Proof.} The decomposition (A1.5) is obtained by splitting the sum in (A1.4) into four regions.\\

\noi {\bf Region 1.} That region is defined by the condition $k \leq j-3$, so that 
$$2^{j-2} \leq |\xi | - |\eta | \leq |\xi - \eta | \leq |\xi | + |\eta | \leq 2^{j+2}$$

\noi and therefore $|j- \ell | \leq 2$. We obtain
$$M_1 = \sum_j \int d\xi \ d\eta (P_1 \widehat{u}_j)(\xi ) \ ( \widehat{S_{j-3} (P_2v)})(\eta ) \ \widehat{\widetilde{m}_j^{(2)}} (\xi - \eta ) \left ( |\xi |^\lambda - |\eta|^\lambda \right ) \ . \eqno({\rm A1.9})$$

\noi The contribution of $|\xi |^\lambda$ is estimated by (A1.1) (A1.2) and the homogeneity of $P_1$, $P_2$ as
$$|M_{1\xi}| \ \leq \ \sum_j \parallel \omega^\lambda P_1 u_j \parallel_{r_1} \ \parallel S_{j-3} (P_2 v) \parallel_{r_2} \ \parallel  \widetilde{m}_j^{(2)} \parallel_{r_0}$$
$$\leq \ \parallel \omega^\lambda P_1 \widetilde{\varphi}_0\parallel_1 \ \parallel P_2 \psi_0 \parallel_1 \sum_j 2^{\mu j-3\alpha_2} \parallel u_j \parallel_{r_1} \ \parallel v \parallel_{r_2} \ \parallel  \widetilde{m}_j^{(2)} \parallel_{r_0}$$
$$\leq \ C \parallel u ; \dot{B}_{r_1, 2}^{\theta \mu}\parallel \ \parallel v \parallel_{r_2}  \ \parallel m ; \dot{B}_{r_0, 2}^{(1 - \theta ) \mu}\parallel \ . \eqno({\rm A1.10})$$

\noi The contribution of $|\eta |^\lambda$ is estimated as
$$|M_{1\eta}| \ \leq \ \sum_j \parallel P_1 u_j \parallel_{r_1} \ \parallel \omega^\lambda S_{j-3} (P_2 v) \parallel_{r_2} \ \parallel  \widetilde{m}_j^{(2)} \parallel_{r_0}$$
$$\leq \ \parallel P_1 \widetilde{\varphi}_0\parallel_1 \ \parallel \omega^\lambda P_2 \psi_0 \parallel_1 \sum_j 2^{\mu (j-3) + 3\alpha_1} \parallel u_j \parallel_{r_1} \ \parallel v \parallel_{r_2} \ \parallel  \widetilde{m}_j^{(2)} \parallel_{r_0}$$

\noi and is therefore estimated by the last member of (A1.10). This proves the estimate (A1.6) for $M_1$.\\

\noi {\bf Region 2.} That region is defined by the condition $j \leq k -3$. The estimate of $M_2$ is then obtained from that of $M_1$ by exchanging $P_1$ with $P_2$ and $u$ with $v$. This proves the estimate (A1.7) for $M_2$. \\

The remaining region $|j-k| \leq 2$ is split again into two regions 3 and 4. The important term is $M_4$ for which the commutator produces a cancellation, and the estimate of $M_4$ requires that $\tau \xi + (1 - \tau )\eta$ stays away from zero for $0 \leq \tau \leq 1$. The harmless term $M_3$ takes care of the situation where that condition is not satisfied. \\

\noi {\bf Region 3.} That region is defined by the conditions $|j - k| \leq 2$ and $\ell \geq j - 4$. The first condition implies that $\ell \leq j + 4$ so that $|\ell - j| \leq 4$ and therefore
$$M_3 = \sum_j \int d\xi \ d\eta \ (P_1 \widehat{u}_j)(\xi )\ ( \widehat{P_2 \widetilde{v}_j^{(2)}} ) (\eta ) \ \widehat{\widetilde{m}_j^{(4)}} (\xi - \eta ) \left ( |\xi |^\lambda - |\eta|^\lambda \right ) \ . \eqno({\rm A1.11})$$

\noi We estimate $M_3$ by
$$|M_{3}| \ \leq \  \left ( \parallel \omega^\lambda P_1 \widetilde{\varphi}_0 \parallel_{1} \ \parallel P_2 \widetilde{\varphi}_0^{(3)} \parallel_{1} \ + \ \parallel P_1 \widetilde{\varphi}_0 \parallel_1\ \parallel  \omega^\lambda P_2 \widetilde{\varphi}_0^{(3)}\parallel_{1}\right )$$
$$\times \sum_j 2^{\mu j}  \parallel u_j\parallel_{r_1} \ \parallel \widetilde{v}_j^{(2)}\parallel_{r_2} \ \parallel  \widetilde{m}_j^{(4)} \parallel_{r_0}$$
$$\leq \ C \left ( \left ( \parallel u ; \dot{B}_{r_1, 2}^{\theta \mu}\parallel \ \parallel v ; \dot{B}_{r_2, \infty}^{0}\parallel\right ) \ \wedge \  \left (  \parallel u ; \dot{B}_{r_1, \infty}^{0}\parallel \ \parallel v ; \dot{B}_{r_2, 2}^{\theta  \mu}\parallel \right ) \right ) \parallel m ; \dot{B}_{r_0, 2}^{(1 - \theta ) \mu}\parallel   \eqno({\rm A1.12})$$

\noi from which the estimates (A1.6) and (A1.7) for $M_3$ follow. \\

\noi {\bf Region 4.} That region is defined by the conditions $|j - k| \leq 2$ and $\ell \leq j - 5$, so that 
$$M_4 = \sum_j \int d\xi \ d\eta (P_1 \widehat{u}_j)(\xi )\ ( \widehat{P_2\widetilde{v}_j^{(2)}} )(\eta ) \  ( \widehat{S_{j-5} (m)}  )  (\xi - \eta ) \left ( |\xi |^\lambda - |\eta|^\lambda \right ) \ . \eqno({\rm A1.13})$$

\noi We rewrite
$$|\xi |^\lambda - |\eta|^\lambda = \lambda \int_0^1 d\tau (\xi - \eta )\cdot (\tau \xi + (1 - \tau ) \eta )\ | \tau \xi + (1 - \tau ) \eta |^{\lambda - 2}$$

\noi so that
$$M_4 = \sum_j \int d\xi \ d\eta \ \widehat{f}_j(\xi , \eta ) \ (P_1 \widehat{u}_j)(\xi )\ ( \widehat{P_2\widetilde{v}_j^{(2)}})(\eta ) \  ( \widehat{\nabla S_{j-5} (m)} )  (\xi - \eta ) \eqno({\rm A1.14})$$

\noi where
$$\widehat{f}_j(\xi , \eta ) = - i \lambda \int_0^1 d\tau \ \widehat{\widetilde{\varphi}_j(\xi )} \ \widehat{\widetilde{\varphi}_j^{(3)}}(\eta ) \widehat{\psi}_{j-4} (\xi - \eta ) (\tau \xi + (1 - \tau )\eta )\ |\tau \xi + (1 - \tau )\eta |^{\lambda - 2} \ . \eqno({\rm A1.15})$$

\noi On the support of $\widehat{f}_j$, we have
$$|\xi - \eta | \leq 2^{j-3} \leq |\xi |/2$$

\noi so that
$$|\tau \xi + (1 - \tau )\eta| \geq |\xi | - |\xi - \eta | \geq |\xi|/2$$

\noi and therefore $\widehat{f}_j \in {\cal C}_0^\infty$. We then estimate by Lemma A1.1 and the fact that $\widehat{f}_j$ is homogeneous of degree $\lambda - 1$
$$|M_{4}| \ \leq \  \sum_j  \parallel f_j \parallel_{1} \ \parallel P_1 u_j  \parallel_{r_1} \ \parallel P_2\widetilde{v}_j^{(2)}\parallel_{r_2} \ \parallel  S_{j-5}(\nabla m) \parallel_{r_0}$$
$$\leq \ \parallel f_0 \parallel_{1} \ \parallel P_1 \widetilde{\varphi}_0\parallel_1 \ \parallel P_2 \widetilde{\varphi}_0^{(3)} \parallel_{1} \ \parallel \omega^{1 - \nu} \psi_0 \parallel_1 \ \sum_j 2^{(\mu - \nu )j} \parallel  u_j  \parallel_{r_1} \ \parallel \widetilde{v}_j^{(2)}\parallel_{r_2} \ \parallel \omega^{\nu - 1} \nabla m \parallel_{r_0}\eqno({\rm A1.16})$$

\noi from which the estimate (A1.8) follows. \par \nobreak \hfill $\sq$\par

\noi {\bf Remark A1.1.} The only assumptions needed on $P_1$, $P_2$ beyond homogeneity are the fact that $P_i\varphi_0$, $\omega^\lambda P_i \varphi_0$, $P_i\psi_0$ and $\omega^\lambda P_i \psi_0$ all belong to $L^1$. If the $P_i$'s are smooth outside of the origin, this is obvious for $\varphi_0$, and follows from the dyadic decomposition of $\psi_0$ if $\alpha_i > 0$. The assumptions made in Lemma 2.4 are trivially sufficient. \par

We now derive a slightly more general lemma than Lemma 2.4.\\

\noi {\bf Lemma A1.3.} {\it Let $\lambda > 0$ and let $P_i$, $i = 1,2$, satisfy the assumptions of Lemma 2.4. Then the estimate (\ref{2.5e}) holds with
$$\left \{ \begin{array}{l} \delta (q_1) = \mu + n/2 - \left ( \sigma_0 + \delta (r_0) + \sigma_2 + \delta (r_2)\right )  \\ \\  \delta (q_2) = \mu + n/2 - \left ( \sigma_0 + \delta (r_0) + \sigma_1 + \delta (r_1)\right )  \\ \\  \delta (q_0) = \mu - \nu + n/2 - \left ( \sigma_1 + \delta (r_1) + \sigma_2 + \delta (r_2)\right )  \end{array} \right .\eqno({\rm A1.17})$$

\noi under the conditions $0 \leq \nu \leq 1$, $1 \leq r_i$, $q_i \leq \infty$ and
$$\left \{ \begin{array}{l} \sigma_0 + \left (\sigma_1 \wedge  \sigma_2\right ) \geq \mu  \\ \\ \sigma_1 +  \sigma_2 \geq \mu - \nu  \end{array} \right .\eqno({\rm A1.18})$$

\noi where $\mu = \lambda + \alpha_1 + \alpha_2$ and $\delta (r) = n/2 - n/r$ for $1 \leq r \leq \infty$.} \\

\noi {\bf Proof.} We rewrite the basic estimates (A1.6)-(A1.8) with a slightly different notation, namely
$$|M_1| + |M_3| \ \leq \ C\parallel u; \dot{B}_{s_1, 2}^{\mu_1} \parallel \ \parallel v \parallel_{q_2} \ \parallel m ; \dot{B}_{s_0, 2}^{\mu_0} \parallel\eqno({\rm A1.19})$$
$$|M_2| + |M_3| \ \leq \ C\parallel u \parallel_{q_1} \ \parallel v ; \dot{B}_{s_2, 2}^{\mu_2}\parallel \ \parallel m ; \dot{B}_{s'_0, 2}^{\mu '_0} \parallel\eqno({\rm A1.20})$$
$$|M_4| \ \leq \ C\parallel u; \dot{B}_{s'_1, 2}^{\mu '_1} \parallel \ \parallel v ; \dot{B}_{s'_2, 2}^{\mu '_2}\parallel \ \parallel \omega^{\nu - 1} \nabla m \parallel_{q_0} \eqno({\rm A1.21})$$

\noi under the conditions $1 \leq s_i, s'_i, q_i \leq \infty$,
$$\left \{ \begin{array}{l} \mu_1 + \mu_0 =  \mu_2 + \mu '_0 = \mu \\ \\ \mu '_1 + \mu '_2 = \mu - \nu  \end{array} \right .\eqno({\rm A1.22})$$

\noi and the H\"older condition
$$1/s_1 + 1/q_2 + 1/s_0 = 1/q_1 + 1/s_2 + 1/s'_0 = 1/s'_1 + 1/s'_2 + 1/q_0 = 1 \ . \eqno({\rm A1.23})$$

\noi We want to choose the parameters $\mu_i$, $\mu '_i$, $s_i$, $s'_i$ in such a way that the $\dot{B}$ norms in (A1.19)-(A1.21) are controlled by the corresponding norms $\dot{B}$ in (\ref{2.5e}) through Sobolev inequalities. This holds provided
$$r_i \leq s_i, s'_i \leq \infty \quad , \qquad i = 0,1,2 \ , \eqno({\rm A1.24})$$

$$\left \{ \begin{array}{l} \mu_1 + \delta (s_1) = \mu{'}_1 +  \delta (s{'}_1)=  \sigma_1 + \delta (r_1) \\ \\  \mu_2 + \delta (s_2) = \mu{'}_2 +  \delta (s{'}_2)=  \sigma_2 + \delta (r_2) \\ \\\mu_0 + \delta (s_0) = \mu{'}_0 +  \delta (s{'}_0)=  \sigma_0 + \delta (r_0) \ . \\ \\ \end{array} \right .\eqno({\rm A1.25})$$

\noi Eliminating $\mu_i$, $\mu '_i$ between (A1.22) and (A1.25) yields
$$\left \{ \begin{array}{l} \sigma_0 + \sigma_1 = \mu +  \delta (s_0) - \delta (r_0) +  \delta (s_1) - \delta (r_1)\\ \\  \sigma_0 + \sigma_2 = \mu +  \delta (s'_0) - \delta (r_0) +  \delta (s_2) - \delta (r_2)\\ \\ \sigma_1 + \sigma_2 = \mu - \nu +  \delta (s'_1) - \delta (r_1) +  \delta (s'_2) - \delta (r_2) \ .  \end{array} \right .\eqno({\rm A1.26})$$

\noi The conditions (A1.24) are then equivalent to
$$\left \{ \begin{array}{l} \mu \leq \sigma_0 + \sigma_1 \leq \mu +  n - \delta (r_0) - \delta (r_1) \\ \\  \mu \leq \sigma_0 + \sigma_2 \leq \mu + n -  \delta (r_0) - \delta (r_2) \\ \\ \mu - \nu \leq \sigma_1 + \sigma_2 \leq \mu - \nu + n - \delta (r_1) - \delta (r_2)  \ .  \end{array} \right .\eqno({\rm A1.27})$$

\noi On the other hand, the H\"older condition (A1.23) can be rewritten as
$$\delta (s_1) + \delta (q_2) + \delta (s_0) = \delta (q_1) + \delta (s_2) + \delta (s'_0) = \delta (s'_1) + \delta (s'_2) + \delta (q_0) = n/2$$

\noi and reduces to (A1.17) by the use of (A1.26). The left hand conditions of (A1.27) coincide with (A1.18) while the right hand conditions reduce to the already imposed conditions $q_i \geq 1$. \par\nobreak \hfill $\sq$ \par

Lemma 2.4 is the special case of Lemma A1.3 where one imposes in addition the global homogeneity condition (\ref{2.7e}) under which (A1.17) reduces to (\ref{2.6e}).

\mysection*{Appendix A2}
\hspace*{\parindent}
\noi {\bf Proof of Proposition 3.1.}

The proof proceeds by a parabolic regularization and a limiting procedure. We consider separately the cases $t \geq t_0$ and $t \leq t_0$ and we begin with $t \geq t_0$. We replace (\ref{1.23e}) by
$$i \partial_t v' = - (1/2) (1 - i \eta ) \Delta v' + \widetilde{L} v' \eqno({\rm A2.1})$$

\noi where $\widetilde{L}$ is defined by (\ref{3.16e}) and $0 < \eta \leq 1$. We recast the Cauchy problem for (A2.1) with initial data $v'(t_0) = v'_0$ in the form of the integral equation
$$v'(t) = U_\eta (t-t_0) v'_0 - i \int_{t_0}^t dt'\ U_\eta (t-t')  \widetilde{L}v'(t')\eqno({\rm A2.2})$$

\noi where
$$U_\eta (t) = \exp \left ( i (t/2) (1 - i \eta ) \Delta \right ) \ . \eqno({\rm A2.3})$$

We first solve (A2.2) locally in time by contraction in ${\cal C}([t_0, t_0 + T_0], H^{\rho '})$ for some $T_0 > 0$. The semi group $U_\eta$ satisfies the estimate
$$\parallel U_\eta (t) \nabla v \parallel \ \leq \ (\eta t)^{-1/2} \parallel v \parallel$$

\noi so that by Lemmas 2.1 and 2.2
\begin{eqnarray*}
&&\parallel \omega^\sigma U_\eta (t-t') \widetilde{L} v'(t') \parallel \\
&&\leq \ (\eta (t - t'))^{-1/2} \parallel \omega^\sigma sv' (t') \parallel \ + \ \parallel \omega^\sigma (\nabla \cdot s) v' (t') \parallel \ + \ \parallel \omega^\sigma fv' (t') \parallel \\
&&\leq \ C \left \{ (\eta (t - t'))^{-1/2} \parallel \omega^{n/2 \pm 0} s \parallel \ + \ \parallel \omega^{n/2 \pm 0} \nabla s \parallel \ + \ \parallel \omega^{n/2\pm 0} f \parallel \right \} \parallel \omega^\sigma v' (t') \parallel 
\end{eqnarray*}
$$\eqno({\rm A2.4})$$

\noi for $0 \leq \sigma \leq \rho '$. We estimate the various terms in the same way as in the proof of Proposition 3.2 (see especially (\ref{3.8e}) (\ref{3.9e})), except for the terms containing $g(v)$ or $g(v_0)$ for which we use the elementary estimate 
$$\parallel \omega^{n/2 \pm 0} g (v_{(0)}) \parallel\ \leq \ C\parallel v_{(0)} ; H^\rho  \parallel^2 \eqno({\rm A2.5})$$

\noi with $\rho > \gamma /2$ instead of the more elaborate estimates (\ref{3.10e})-(\ref{3.12e}). We can then continue (A2.4) as
$$\cdots \leq \ C \Big\{ (\eta (t - t'))^{-1/2} a_0^2 \ t{'}^{\lambda_0 - 1} + a_0^2 \left ( t{'}^{\lambda_1 - 1} + a_0^2\ t{'}^{2\lambda_0 - 2}\right )$$ $$ + \ \left ( a_0^2 \ + \ \parallel v(t') ; H^\rho  \parallel^2 \right ) t^{\gamma - 2} \Big \} \parallel \omega^\sigma v'(t') \parallel \eqno({\rm A2.6})$$

\noi where $a_0 = \parallel v_0 ; H^\rho  \parallel$. It then follows from (A2.6) that (A2.2) can be solved by contraction in ${\cal C} ([t_0, t_0 + T_0], H^{\rho '})$ for $T_0$ sufficiently small. By a standard argument using the linearity of (A2.2), one can extend the solution to $[t_0, T]$. Let $v'_\eta$ be that solution. \par

We next take the limit where $\eta$ tends to zero. For that purpose we first estimate $v'_\eta$ in $L^\infty ([t_0, T], H^{\rho '})$ uniformly in $\eta$. In the same way as in the proof of Proposition 3.2, we estimate
$$\partial_t \parallel \omega^\sigma v'_\eta (t) \parallel^2 \ = \ - \eta \parallel \omega^\sigma \nabla v'_\eta  (t) \parallel^2 \ + \ {\rm Im} \ <v'_\eta (t) , [\omega^{2\sigma} , L(v)] v'_\eta (t)>$$
$$\leq \ C \left ( \parallel \omega^{n/2} \nabla s \parallel \ + \ \parallel \nabla s \parallel_\infty \ + \ \parallel \omega^{n/2} f \parallel \right ) \parallel \omega^\sigma v'_\eta (t) \parallel^2$$  
$$\leq \ N_0 (t) \parallel \omega^\sigma v'_\eta (t) \parallel^2 \eqno({\rm A2.7})$$

\noi with 
$$N_0 (t) = C \left \{ a_0^2 \left ( t^{\lambda_1 - 1} + a_0^2\ t^{2\lambda_0 - 2}\right ) + \left ( a_0^2 + \ \parallel v(t) ; H^\rho  \parallel^2 \right ) t^{\gamma - 2}\right \} \eqno({\rm A2.8})$$

\noi for $0 \leq \sigma \leq \rho '$. Here we have again used (A2.5). It follows from (A2.7) by integration that $v'_\eta$ is estimated in $L^\infty ([t_0, T], H^{\rho '})$ uniformly in $\eta$. By compactness, one can find a sequence of $v'_\eta$ for $\eta$ tending to zero which converges in the weak-$\star$ sense to a limit $v' \in L^\infty ([t_0, T], H^{\rho '})$. The limit $v'$ satisfies the equation (\ref{1.23e}), so that $v' \in {\cal C} ([t_0, T], H^{\rho '- 2}) \cap {\cal C}_w ([t_0, T], H^{\rho '})$. Furthermore $v'_\eta$ tends to $v'$ weakly in $H^{\rho '}$ pointwise in $t$ so that $v'(t_0) = v'_0$. A similar argument yields the same results for $t \leq t_0$. Uniqueness follows from $L^2$-norm conservation and linearity. Strong continuity in $H^{\rho '}$ follows from the estimate (A2.7) which implies continuity of $\parallel v'(t) ; H^{\rho '}\parallel$ at $t_0$ and from uniqueness through a change of initial time.\par\nobreak \hfill $\sq$ \par

\end{document}